\documentclass[12pt]{elsarticle}

\makeatletter
\def\ps@pprintTitle{%
	\let\@oddhead\@empty
	\let\@evenhead\@empty
	\let\@oddfoot\@empty
	\let\@evenfoot\@oddfoot
}
\makeatother

\usepackage[a4paper, total={6.1in, 10in}]{geometry}
\usepackage{amsmath,amssymb,amsthm}
\usepackage{color}
\usepackage{graphicx}
\usepackage{verbatim}
\usepackage{physics}
\usepackage{algorithm}
\usepackage{algpseudocode}
\usepackage{datetime}
\usepackage{epstopdf}

\newtheorem{theorem}{Theorem}[section]

\newtheorem{definition}[theorem]{Definition}

\newtheorem{remark}[theorem]{Remark}

\newcommand{\C}{\mathbb C}

\newcommand{\etal}{\mbox{\emph{et al.\ }}}

\newcommand{\OcD}{{}^c\!D}

\newcommand{\cD}{{}^c\!\partial}

\newcommand{\wfun}[2]{{}_{#1}\Psi_{#2}}

\begin{document}

\begin{frontmatter}
\title{
Highly Accurate Global Pad\'e Approximations of Generalized 
Mittag-Leffler Function and its Inverse
}

\author[kfupm,isarumi]{Ibrahim O. Sarumi}
\author[kfupm,kmfurati]{Khaled M. Furati}
\author[mtsu]{Abdul Q. M. Khaliq}

\address[kfupm]{King Fahd University of Petroleum \& Minerals \\ 
Department of Mathematics \& Statistics \\ 
Dhahran 31261, Saudi Arabia}

\address[isarumi]{isarumi@kfupm.edu.sa}
\address[kmfurati]{kmfurati@kfupm.edu.sa}

\address[mtsu]{Middle Tennessee State University \\ 
Department of Mathematical Sciences\\ 
Murfreesboro, TN 37132-0001, USA \\ 
Abdul.Khaliq@mtsu.edu}

\date{\today, \currenttime}



\begin{abstract}

The two-parametric Mittag-Leffler function (MLF), $E_{\alpha,\beta}$, is fundamental to the study and simulation of fractional differential and integral equations. 
However, these functions are computationally expensive and their numerical 
implementations are challenging.  
In this paper, we present a unified framework for developing global rational approximants of 
$E_{\alpha,\beta}(-x)$, $x>0$, with 
$\{ (\alpha,\beta): 0 < \alpha \leq 1, \beta \geq \alpha, (\alpha,\beta) \ne (1,1) \}$.
This framework is based on the series definition and the asymptotic expansion at infinity.
In particular, we develop three types of fourth-order global rational approximations and discuss how they could be used to approximate the inverse function.
Unlike existing approximations which are either limited to MLF of one parameter 
or of low accuracy for the two-parametric MLF, our rational approximants are of fourth order accuracy 
and have low percentage error globally. 
For efficient utilization, we study the partial fraction decomposition and 
use them to approximate the two-parametric MLF with a matrix argument which arise in the solutions 
of fractional evolution differential and integral equations.

\end{abstract}
	
\end{frontmatter}

\section{Introduction}
In this paper, we consider the two-parametric MLF
\begin{equation}
\label{eq:mlf2}
E_{\alpha, \beta}(z) = \sum_{k=0}^{\infty} \frac{z^k}{\Gamma(\alpha k + \beta)}, 
\qquad \Re \alpha > 0, \quad \beta \in \C, \quad z \in \C.
\end{equation}
This entire function generalizes the MLF of one-parameter, $E_\alpha = E_{\alpha, 1}$.

The function $E_{\alpha, \beta}$ plays a key role in the study and simulation 
of history-dependent evolution models that arise in many engineering and science areas 
such as flow in porous media, pattern recognition, rheology, anomalous diffusion, 
electric networks, etc.
In particular, it is the cornerstone of the development of generalized exponential time differencing (GETD) schemes \cite{Garrappa:2011} 
which extend the notion of exponential integrator \cite{Minchev:2005} to 
time-fractional problems. 

Evaluation of $E_{\alpha,\beta}$ with scalar arguments is very expensive and challenging.
Although the series \eqref{eq:mlf2} converges analytically for all $z\in \C$, 
it is not practical or may not be valid to use it computationally for $|z|\geq 1$. 
Consequently, different techniques for the evaluation of $E_{\alpha,\beta}$ have been developed.
%
Gorenflo, Loutchko and Luchko \cite{Gorenflo:2002} 
proposed an algorithm based on using appropriate techniques for different regions of $\C$.
For small and large $\abs{z}$ values, they used the series definition \eqref{eq:mlf2} 
and the asymptotic series at infinity, respectively. For the intermediate regions they used the integral representations.
A similar approach has been followed by Hilfer and Seybold \cite{Hilfer:2006}.
Garrappa \cite{Garrappa:2015} provided an approach based on the numerical inversion of the Laplace transform. 
For efficient implementation, he 
provided an algorithm for finding the optimal parabolic contour on the basis of the distance and strength 
of the singularities of the Laplace transform. 

The evaluation of MLF with a matrix argument is still a tricky and tough task.
Garrappa \cite{Garrappa:2018} developed an algorithm based on the similarity transform.
This approach requires the evaluation of MLF and its derivatives for each eigenvalue, which is again obtained using the Laplace 
transform. Clearly, massive calculations will be required for large full matrices.

In summary, all existing algorithms for evaluating $E_{\alpha,\beta}$ suffer from some drawbacks such as nontrivial software 
implementation, long CPU time especially when a fine error tolerance is imposed, 
overflow numbers, and catastrophic cancellations. 
Due to these computational complexities and the need for efficient matrix function evaluation, accurate and efficient approximations are imperative.

To the best of authors' knowledge, there have been few studies about rational approximations of MLF. 
Freed \etal \cite{Freed:2002} developed a piecewise approximant for $E_\alpha(-x^\alpha)$, $x>0$, 
based on the truncated series representation for small values, the asymptotic 
expansion for large values, and a Pad\'e type approximant for the intermediate values.
For $E_\alpha(-x)$, $x>0$, 
Starovoitov and Starovoitova \cite{Starovoitov:2007} analyzed Pad\'e type approximants of the form 
$p_n/q_m$, $m \le n$, and discussed their asymptotic rate of convergence on the compact unit disk as 
$n\to\infty$.
Borhanifar and Valizadeh \cite{Borhanifar:2015} constructed a fourth order Pad\'e approximant and used it to develop a numerical scheme 
for the time-space diffusion equation.
Iyiola \etal \cite{Iyiola:2018real} constructed a second order 
non-Pad\'e type rational approximation for $E_{\alpha, \beta}$ using real distinct poles (RDP).
However, although the approximants in \cite{Borhanifar:2015} and \cite{Iyiola:2018real} might be adequate for small values, they fail to 
account for the asymptotic power law behavior.

The global Pad\'e approximation technique introduced by Winitzki \cite{Winitzki:2003} 
has been applied recently to construct rational approximations for the Mittag-Leffler 
function and its generalization.  
In this technique, rational approximations are constructed by matching them with selected
combinations of the series definition and the asymptotic expansion.
Atkinson and Osseiran \cite{Atkinson:2011} used this technique to construct a second-order rational 
approximation for $E_{\alpha}$. 
Later, Ingo \etal \cite{Ingo:2017} showed that the rational approximant in \cite{Atkinson:2011} is not satisfactory for $\alpha$ close to one.
Alternatively, they constructed a fourth-order global approximation for $E_{\alpha}$ that behaves reasonably well for all values of $\alpha \in (0,1)$.
Zeng \cite{Zeng2015} extended this technique to construct a second-order global Pad\'e approximant for $E_{\alpha,\beta}$.
However, this approximation is not satisfactory for $\alpha$ close to one, especially when $\beta = 1$ and it is malfunctioning for  $\beta = \alpha$, $0.5 \leq \alpha < 1$.

As for the inverse of MLF,
Hilfer and Seybold \cite{Hilfer:2006} introduced the inverse of $E_{\alpha,\beta}(z)$ as the solution of the equation 
\begin{equation}
\label{eq:inverse mlf}
L_{\alpha,\beta}(E_{\alpha,\beta}(z)) = z, \qquad z\in \C,
\end{equation}
where $L_{\alpha,\beta}(z)$ is evaluated by solving the functional equation \eqref{eq:inverse mlf} numerically.
They discussed the principal branch of the function and showed that it reduces to the 
principal branch of the logarithm function as $\alpha \to 1$ when $\beta = 1$. 
Hanneken and Achar \cite{Hanneken:2014} proposed a finite series representation for $L_{\alpha, \beta}$ but only for some values of $\alpha \in (0,\frac{1}{2})$ and $\beta = 1, 2$.
Lately, approximations of $L_{\alpha, \beta}$ have been introduced to overcome the difficulty of solving the functional equation~\eqref{eq:inverse mlf}. 
Atkinson and Osseiran \cite{Atkinson:2011}; and Ingo \etal 
\cite{Ingo:2017} discussed the approximation of $(-L_{\alpha,1})$ based on the inversion of their global Pad\'e approximants. 
Similarly, Zeng and Chen \cite{Zeng2015}; and Iyiola \etal \cite{Iyiola:2018real} inverted their second order 
approximants of $E_{\alpha, \beta}$ to obtain an approximation of $(-L_{\alpha, \beta})$.

Consequently, based on the current state of the literature, more accurate rational approximations of $E_{\alpha,\beta}$ and its inverse are needed.
Such approximations are expected to ease computation cost and yield accurate values globally. 
In this paper, we introduce a framework that unifies the notion of global Pad\'e approximation 
for the two-parametric MLF. 
Moreover, we develop different types of fourth-order global rational approximations 
for $E_{\alpha, \beta}(-x)$, $x>0$, for
 $\{ (\alpha,\beta): 0 < \alpha \leq 1, \beta \geq \alpha, (\alpha,\beta) \ne (1,1) \}$.
We also discuss analytically and numerically the approximation errors. 
Furthermore, we present the partial fraction decomposition of the rational approximants together with its advantage in 
efficient implementation for matrix arguments. An algorithm for computing $(-L_{\alpha, \beta})$ based on the inversion 
of our fourth order approximants is presented. All along, we demonstrate through numerical experiments and comparisons, 
that the new developed fourth-order approximants provide superior global approximations for $E_{\alpha,\beta}$ and 
its special case $E_\alpha$.
	
This paper is organized as follows.
Section \ref{sec:GlobalPade} contains the unified framework for the 
global Pad\'e approximation and the error analysis.
In section \ref{sec:second-order GPA}, we discuss the second order global Pad\'e approximant constructed 
by Zeng and Chen in \cite{Zeng2015} and the need for more accurate approximants. 
Section \ref{sec:fourth-order GPA} contains the construction of our fourth order approximants. 
The partial fraction decomposition and the algorithms to 
compute the poles and weights are discussed in section \ref{sec:partial-frac}.
In section \ref{sec:Inv-mlf} we discuss the 
inverse MLF and its approximation through the inversion of our rational 
approximants. 
Graphical and numerical demonstrations of the performance of our approximants are presented in section~\ref{sec:comparison}.
In section \ref{sec:applcications}, we apply our approximants to the solutions of fractional differential and integral equations and systems.

The computations in this paper are performed using Matlab software on a dell laptop with a core i5 processor. 


\section{Global rational approximation for $E_{\alpha,\beta}(-x)$, $x>0$}
\label{sec:GlobalPade}

In this section, we introduce and outline the construction of the global Pad{\'e} approximation 
for $E_{\alpha, \beta}(-x)$, $x>0$, for the cases
\begin{equation}
\label{eq:alpha-beta}
\mathcal{A} = \{ (\alpha,\beta): 0 < \alpha \leq 1, \beta \geq \alpha, (\alpha,\beta) \ne (1,1) \}.
\end{equation}
Our approach is based on the technique proposed by Winitzki in \cite{Winitzki:2003}.
This technique relies on the asymptotic expansion given by the following 
theorem (\cite{Podlubny:1999a}, Theorem 1.4).

\begin{theorem}
Let $\alpha \in (0,2)$, $\beta \in \mathbb{C}$ and $\mu \in \mathbb{R}$,  
$\frac{\pi \alpha}{2} < \mu <\min\{\pi, \pi\alpha\}$. 
Then for $\mu \le \abs{\arg z} \le \pi$, 
\begin{equation}
\label{eq:mlf asymptotic}
E_{\alpha, \beta}(z) =  - \sum_{k=1}^{n}\dfrac{(z)^{-k}}
{\Gamma(\beta - \alpha k)} + \mathcal{O}(|z|^{-(n+1)}),  
\quad \mbox{as } |z| \to \infty, \quad n \geq 1.
\end{equation}
\end{theorem}
In particular, when $\beta = \alpha$, the series in \eqref{eq:mlf asymptotic} takes the form
\begin{align}
\label{eq:mlf asymptotic aa}
E_{\alpha, \alpha}(z) = - \sum_{k=1}^{n-1} \dfrac{ z^{-(k+1)}}{\Gamma(- \alpha k )} 
+ \mathcal{O}(|z|^{-(n+1)}),
\quad \text{as } |z| \to \infty, \quad n \geq 2.
\end{align}
As an abbreviation for the rest of the paper, the cases $\beta > \alpha$ and $\beta = \alpha$ are to be understood as sub-cases of \eqref{eq:alpha-beta}.

\subsection{Definition}
We proceed by considering the function
\begin{equation}
\label{eq:scritptE}
\mathcal{E}_{\alpha,\beta}(x) = s_{\alpha, \beta}(x)E_{\alpha,\beta} (-x),
\end{equation}
with $s_{\alpha, \beta}(x)$ chosen so that the first term in the asymptotic expansion of 
$\mathcal{E}_{\alpha,\beta}$ is 1. 
It follows from \eqref{eq:mlf asymptotic} and 
\eqref{eq:mlf asymptotic aa} that
\begin{equation}
\label{eq:s_ab}
s_{\alpha, \beta}(x) = 
\begin{cases}
\Gamma(\beta - \alpha) x,  & \beta > \alpha, 
\\ &
\\
-\Gamma( - \alpha) x^2,    & \beta = \alpha.
\end{cases}
\end{equation}
This function admits the following behavior:
\begin{equation}
\label{eq:global-pade}
\mathcal{E}_{\alpha,\beta}(x) = 
\begin{cases}
a(x) + \mathcal{O}(x^m),  & \text{ as } x \to 0, \quad m \geq 
\begin{cases} 2, & \beta > \alpha,
\\
3, &  \beta = \alpha,
\end{cases}
\\ \\
b(x^{-1}) + \mathcal{O}(x^{-n}),  & \text{ as } x \to \infty, \quad n \ge
\begin{cases} 1, & \beta > \alpha,
\\
2, &  \beta = \alpha,
\end{cases}
\end{cases}
\end{equation}
where, from \eqref{eq:mlf2},
\begin{equation}
a(x) =
\begin{cases} \displaystyle
\Gamma(\beta - \alpha)\,x \, \sum\limits_{k=0}^{m-2}
\frac{(-x)^k}{\Gamma(\beta + \alpha k)}, \qquad &  \beta > \alpha,
\\ \\ \displaystyle
-\Gamma(-\alpha)\,x^2 \, \sum\limits_{k=0}^{m-3}
\frac{(-x)^k}{\Gamma(\alpha k + \alpha)}, &  \beta = \alpha,
\end{cases}
\end{equation} 
and from \eqref{eq:mlf asymptotic}--\eqref{eq:mlf asymptotic aa},
\begin{equation}
b(x^{-1}) = 
\begin{cases}
\displaystyle
-\Gamma(\beta - \alpha)\,x\sum\limits_{k=1}^{n}
\frac{(-x)^{-k}}{\Gamma(\beta - \alpha k)}, \qquad & \beta > \alpha,
\\ \\ \displaystyle
\Gamma( - \alpha) \, x^2\sum\limits_{k=1}^{n}
\frac{(-x)^{-(k+1)}}{\Gamma( - \alpha k)}, &  \beta = \alpha.
\end{cases}
\end{equation}

Note that when $n=1$, then $b(x^{-1}) = 1$. 
We will see later (equation \eqref{eq:asympt-coeff}) that in this case the asymptotic expansion 
in \eqref{eq:global-pade} does not contribute to the rational approximation 
of $\mathcal{E}_{\alpha,\beta}$.
Therefore, for our purposes, we always take $n>1$.

Next, we introduce the following definition.

\begin{definition}
\label{def:rmn} 
Consider $E_{\alpha,\beta}$ with $(\alpha,\beta) \in \mathcal{A}$.
Let $m$ and $n$ be positive integers such that 
\begin{equation}
\label{eq:mn-values}
n > 1, \qquad 
m \geq 
\begin{cases} 
	2, & \text{ if } \quad \beta > \alpha,
	\\
	3, &  \text{ if } \quad \beta = \alpha,
\end{cases}, \qquad m+n \text{ is odd}.
\end{equation}
Then, the global Pad{\'e} approximation, $R_{\alpha, \beta}^{m, n}(x)$, 
of type $(m,n)$ for $E_{\alpha,\beta} (-x)$ 
is defined as
\begin{equation}
\label{eq:rmn}
R_{\alpha, \beta}^{m, n}(x) = \dfrac{1}{s_{\alpha,\beta}(x)}\dfrac{p(x)}{q(x)}, \quad 
 0 < \alpha \leq 1, \quad \beta \geq \alpha,  \quad (\alpha, \beta) \ne (1,1),
\end{equation}
where $p$ and $q$ are polynomials of degree $\nu$, 
\begin{equation}\label{eq:nu}
\nu := \frac{m + n - 1}{2} \geq 1,
\end{equation} 
such that $q(x) \ne 0$ for $x>0$ and 
\begin{equation}
\label{eq:pq condition}
\dfrac{p(x)}{q(x)} = 
\begin{cases}
a(x) + \mathcal{O}(x^{m-\nu}),  & \text{ as } x \to 0, 
\\ &
\\
b(x^{-1}) + \mathcal{O}(x^{-n}),    & \text{ as } x \to \infty.
\end{cases}
\end{equation}
\end{definition}

Next, we present the procedure for constructing $R^{m,n}_{\alpha,\beta}(x)$.

\subsection{Construction of $R^{m,n}_{\alpha, \beta}(x)$}
\label{subsec:gpa construction}

Let $m$ and $n$ be as in \eqref{eq:mn-values}. We seek a rational approximation of the form
\begin{equation}
\label{eq:gpa general form}
\mathcal{E}_{\alpha,\beta}(x)  \approx  \dfrac{p(x)}{q(x)} =
\dfrac{p_0 + p_1 x + \dots +p_\nu x^\nu}{q_0 + q_1 x + \dots + q_\nu x^\nu},
\end{equation} 
where $\nu$ is as in \eqref{eq:nu}.
This means that $2\nu+1$ coefficients are to be determined.

Since $\mathcal{E}_{\alpha, \beta}(x) \to 1$ as $x \to \infty$ and 
$\lim_{x\to \infty} p(x)/q(x) = p_\nu/q_\nu$, we can set 
$$
p_\nu = q_\nu = 1.
$$
To find the other $2\nu$ unknowns $\{p_i,\, q_i\}_{i=0}^{\nu-1}$, we solve the system
of linear equations obtained by satisfying the requirement \eqref{eq:pq condition} 
which takes the form
\begin{align}
\label{eq:gpa matching}
p(x) - q(x)a(x) & = \mathcal{O}(x^m), & \text{ as } x \to 0,
\\ 
\label{eq:gpa matching2}
x^{-\nu}p(x) - x^{-\nu}q(x)b(x^{-1}) & = \mathcal{O}(x^{-n}), & \text{ as } x \to \infty.
\end{align}

By expanding the left-hand side of \eqref{eq:gpa matching}, it follows that the coefficients of 
$x^k$, $k = 0, 1, \ldots, m-1$, must vanish. 
As such, we obtain $m$ linear equations. 
Similarly, by expanding the left hand side of \eqref{eq:gpa matching2}, the coefficients of
\begin{equation}
\label{eq:asympt-coeff}
x^{-1}, x^{-2}, \dots, x^{-(n-1)}
\end{equation}
must vanish and we obtain another system of $n-1$ linear equations.
Collectively, \eqref{eq:gpa matching} and \eqref{eq:gpa matching2} yield a linear system of $m+n-1$ (= $2\nu$) equations which are then solved for the $2\nu$ unknowns. 
By inspection, we have
\begin{equation}\label{eq:initial p values}
\begin{cases}
p_0 = 0, & \beta > \alpha , 
\\
p_0 = p_1 = 0,  \quad  &  \beta = \alpha.
\end{cases}
\end{equation}
Hence, for $\beta>\alpha$ we solve $2\nu -1$ ($= m+n-2$) equations with $m+n-2$ unknowns,
while for $\beta = \alpha$ we solve $2\nu -2$ (= $m+n-3$) equations with $m+n-3$ unknowns.

We will see later (Remark \ref{remark:error}) that for controlling the approximation error we must have $m > n$.
Table \eqref{tab:mn combinations} provides the order of approximations for the types $(m,n)$ with $m>n > 1$ and $m+n$ is odd. 

\begin{table}[]
\caption{Order of approximations for different types $(m,n)$}
\centering
\label{tab:mn combinations}
\bigskip
\begin{tabular}{c|cccc} 
\hline \\
\qquad $n$ & \quad 2 \quad & \quad 3 \quad & \quad 4 \quad & \quad 5 \quad \\ 
$m$ & & & & \\ \hline \\
3& 2 &   &  &  \\
4& $\cdot$ & 3 &  & \\
5& 3 & $\cdot$ & 4 & \\
6& $\cdot$ & 4 & $\cdot$ & 5\\
7& 4 & $\cdot$ & 5 & $\cdot$\\
8& $\cdot$ & 5 & $\cdot$ & 6 \\ \\\hline
\end{tabular}

\end{table} 

\subsection{Approximation error}

Consider the pointwise error
\begin{equation}
\label{eq:approx error}
e_{\alpha, \beta}^{m,n}(x)  := E_{\alpha, \beta}(-x) - R_{\alpha,\beta}^{m, n}(x), \qquad x>0.
\end{equation} 
Equations \eqref{eq:scritptE}, \eqref{eq:s_ab}, \eqref{eq:global-pade}, \eqref{eq:rmn}, 
and \eqref{eq:pq condition} yield the following orders.
As $x \to 0$, we have
\begin{align}
\nonumber
e_{\alpha, \beta}^{m,n}(x) = E_{\alpha, \beta}(-x) - R_{\alpha,\beta}^{m, n}(x) 
& = 
\frac{1}{s_{\alpha,\beta}(x)} \left\{ \mathcal{E}_{\alpha, \beta}(x)  - \frac{p(x)}{q(x)}
\right\}
\\ &=  \nonumber
\frac{1}{s_{\alpha,\beta}(x)} \left\{ a(x) + \mathcal{O}(x^m)  - a(x) - \mathcal{O}(x^{m-\nu}) \right\}
\\ & =  \nonumber
\frac{1}{s_{\alpha,\beta}(x)} \left\{ \mathcal{O}(x^m)  + \mathcal{O}(x^{m-\nu}) \right\}
\\ \nonumber \\ & = 
\label{eq:error near zero}
\begin{cases}
\mathcal{O} (x^{m-\nu-1}),  \qquad & \beta > \alpha, 
\\ 
\mathcal{O} (x^{m-\nu-2}),    & \beta = \alpha.
\end{cases}
\end{align}
As $x\to\infty$, we have
\begin{align}
\nonumber
e_{\alpha, \beta}^{m,n}(x) = E_{\alpha, \beta}(-x) - R_{\alpha,\beta}^{m, n}(x) 
& = 
\frac{1}{s_{\alpha,\beta}(x)} \left\{ \mathcal{E}_{\alpha, \beta}(x)  - \frac{p(x)}{q(x)}
\right\}
\\ &=  \nonumber
\frac{1}{s_{\alpha,\beta}(x)} \left\{ b(x^{-1}) + \mathcal{O}(x^{-n}) - b(x^{-1}) - \mathcal{O}(x^{-n})
\right\}
\\ & = \nonumber
\frac{1}{s_{\alpha,\beta}(x)} \, \mathcal{O}(x^{-n})
\\ \nonumber \\ & =
\label{eq:error at infinity}
\begin{cases}
\mathcal{O} (x^{-n-1}),  \qquad & n > 1, \quad \beta > \alpha, 
\\ 
\mathcal{O} (x^{-n-2}),    &  n>1, \quad \beta = \alpha.
\end{cases}
\end{align}

\begin{remark}
\label{remark:error}
By inspection of \eqref{eq:error near zero} and \eqref{eq:error at infinity}, one can observe the following.

\begin{itemize}

\item 
For reliable approximations of $R_{\alpha, \beta}^{m,n}$ at small values, one should consider $m \geq n+1$ when $\beta \ne \alpha$ and $m \geq n+3$ when $\beta = \alpha$. 
This is why, for example, $R_{\alpha, \beta}^{5,4}$ is not a good approximation when $\beta = \alpha$.

\item
For large values of $x$, 
the approximation error can be made arbitrary small by taking $n$ sufficiently large.
However, this may not be the case since the asymptotic series 
\eqref{eq:mlf asymptotic} of the Mittag-Leffler function  is divergent. 

\end{itemize} 
\end{remark}

\begin{remark}
In all the numerical experiments and comparisons throughout this paper, the term "exact" values of MLF refers to the values computed using the routines discussed in \cite{Garrappa:2015} and \cite{Garrappa:2018}.
\end{remark}

\section{Second-order global Pad\'e approximant $R^{3,2}_{\alpha,\beta}$}
\label{sec:second-order GPA}
For completeness, we provide here an overview of the second order global Pad\'e approximant  
$R^{3,2}_{\alpha, \beta}(x)$ constructed by Zeng and Chen in \cite{Zeng2015}. The approximant is given by
\begin{equation}
\label{eq:Zeng gpa}
R_{\alpha, \beta}^{3, 2}(x) = \dfrac{1}{\Gamma(\beta - \alpha)} \,
\dfrac{p_1 + x }{q_0 + q_1 x +  x^2 }, \qquad \beta > \alpha, 
\end{equation}
with 
\begin{equation}
\begin{aligned}
p_1  & = c_{\alpha,\beta} \left[
\Gamma(\beta) \Gamma(\beta + \alpha ) - 
\frac{\Gamma(\beta + \alpha)\Gamma^2(\beta - \alpha)}{\Gamma(\beta - 2\alpha)}
\right], 
\\ 
q_0  & = c_{\alpha,\beta} \left[\
\frac{\Gamma^2(\beta) \Gamma(\beta + \alpha )}{\Gamma(\beta - \alpha)}
- \frac{\Gamma(\beta) \Gamma(\beta + \alpha)\Gamma(\beta - \alpha)}
{\Gamma(\beta - 2\alpha)}
\right],
\\ 
q_1 & = c_{\alpha,\beta} \left[
\Gamma(\beta) \Gamma(\beta + \alpha ) - 
\frac{\Gamma^2(\beta)\Gamma(\beta - \alpha)}{\Gamma(\beta - 2\alpha)}
\right], 
\\
c_{\alpha,\beta} &= \frac{1}{\Gamma(\beta + \alpha) \Gamma(\beta - \alpha) - \Gamma^2(\beta)},
\end{aligned}
\end{equation}
and 
\begin{equation}
R_{\alpha, \alpha}^{3, 2}(x) = 
\dfrac{\alpha}{\Gamma(1+\alpha) + 
\frac{2\Gamma(1 - \alpha)^2}{\Gamma(1-2\alpha)} \, x + 
\Gamma(1 - \alpha) \, x^2},
\qquad 0 < \alpha < 1.
\end{equation}

\begin{figure}
\centering
	\includegraphics[width=.49\textwidth]{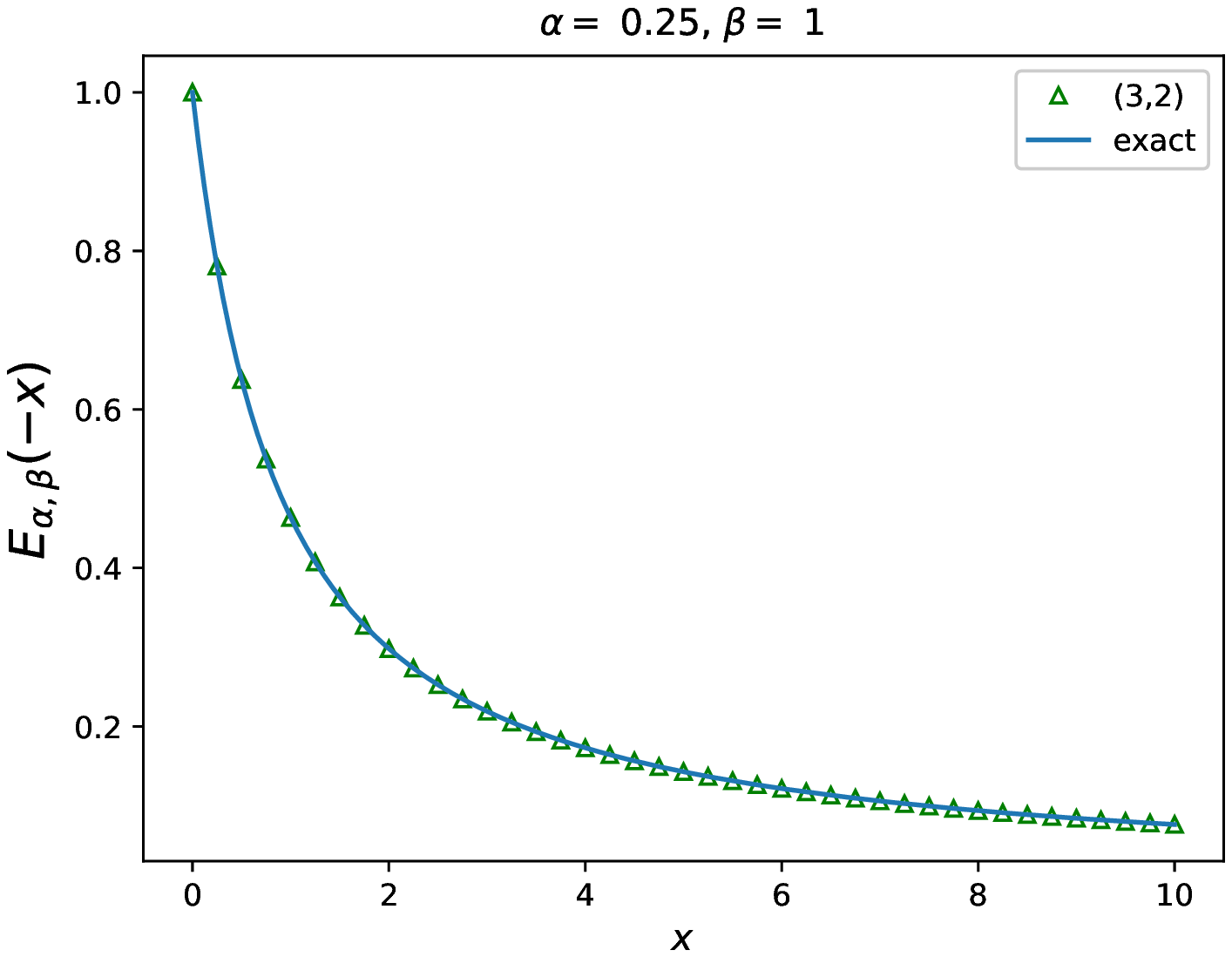}
	\includegraphics[width=.49\textwidth]{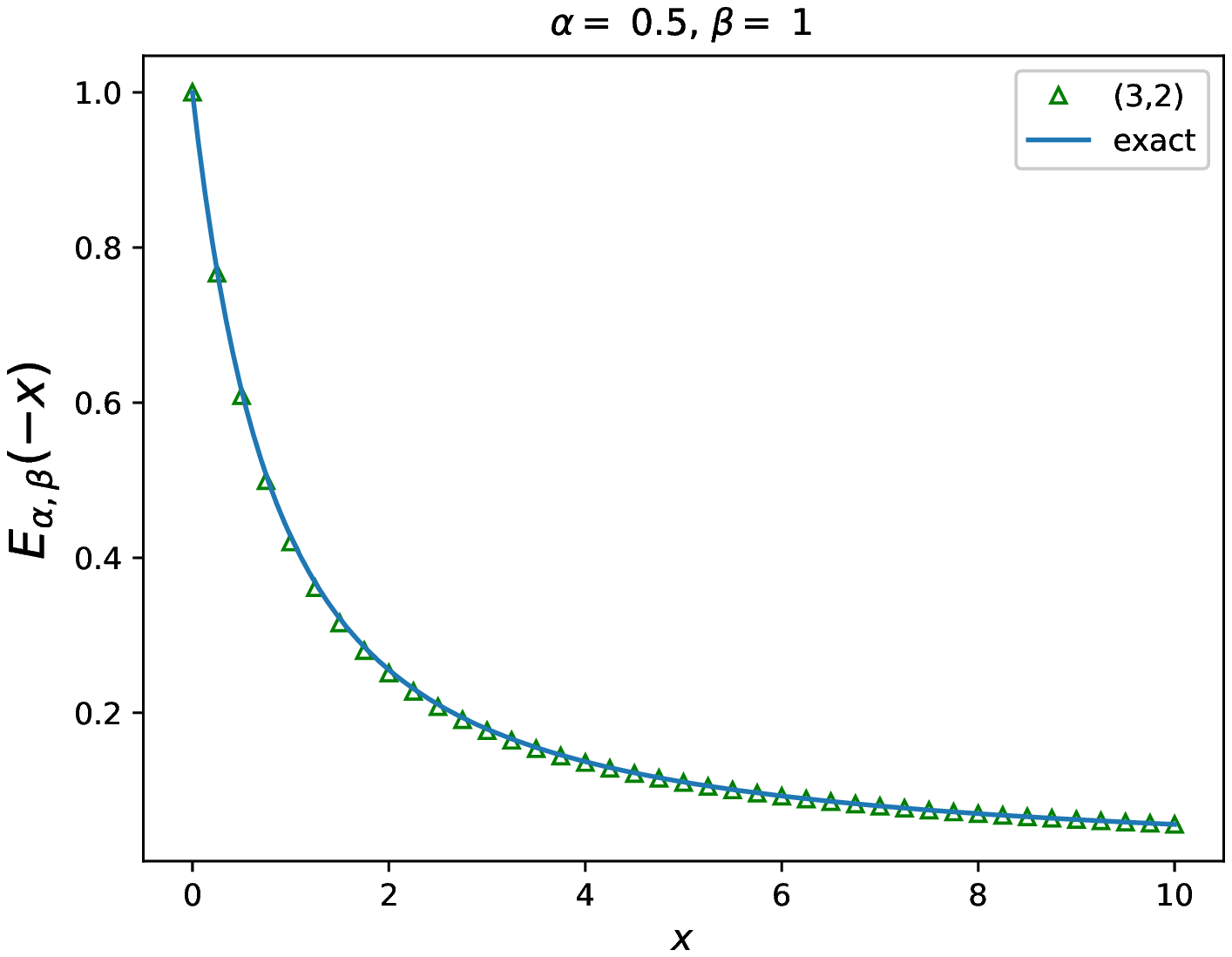}
	\\
	\includegraphics[width=.49\textwidth]{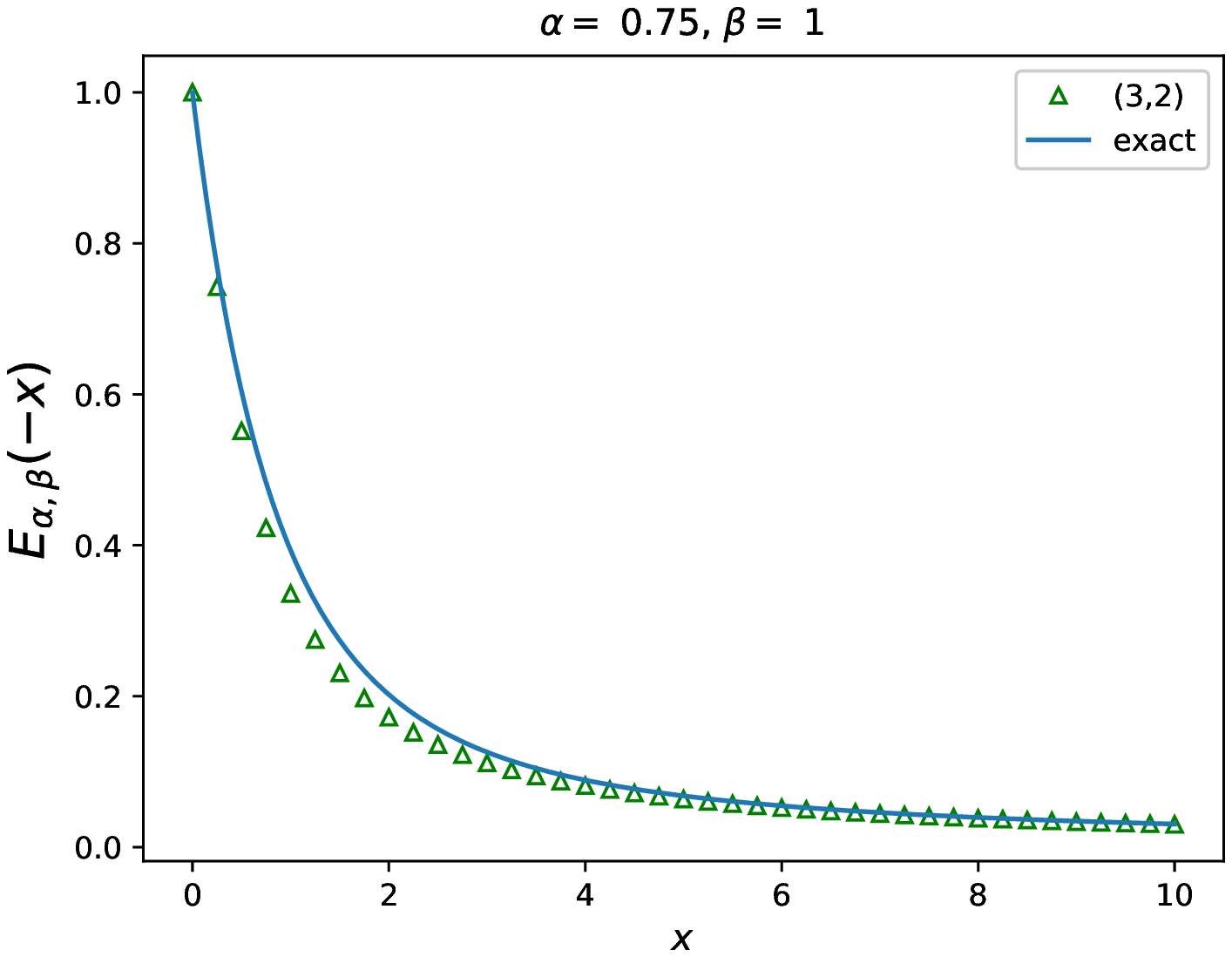}
	\includegraphics[width=.49\textwidth]{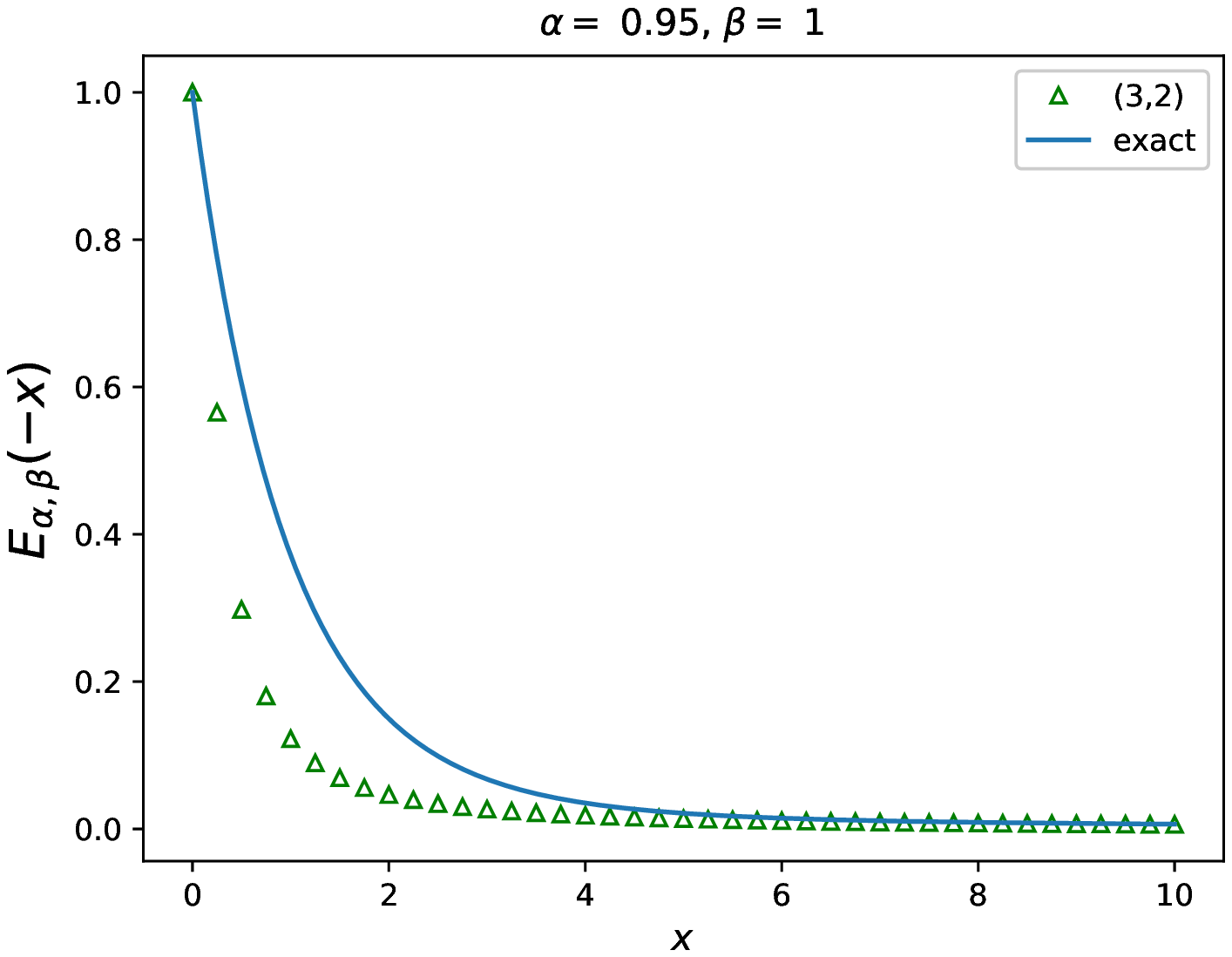}
\caption{Plots of $R^{3,2}_{\alpha,1}$ for different values of $\alpha$}
\label{fig:Zeng small alpha beta 1}
\end{figure}

\begin{figure}
\centering
\includegraphics[width=.49\textwidth]{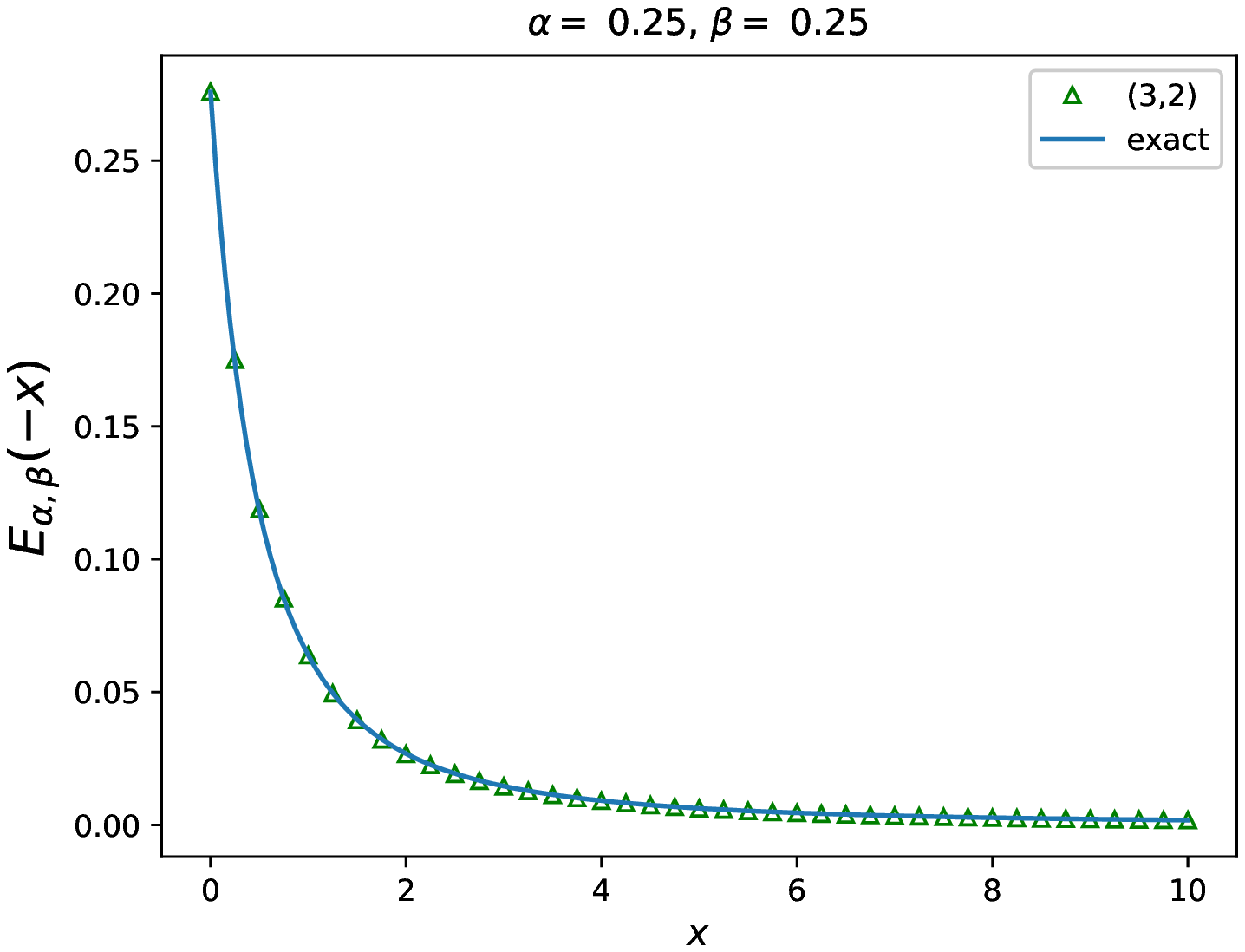}
\includegraphics[width=.49\textwidth]{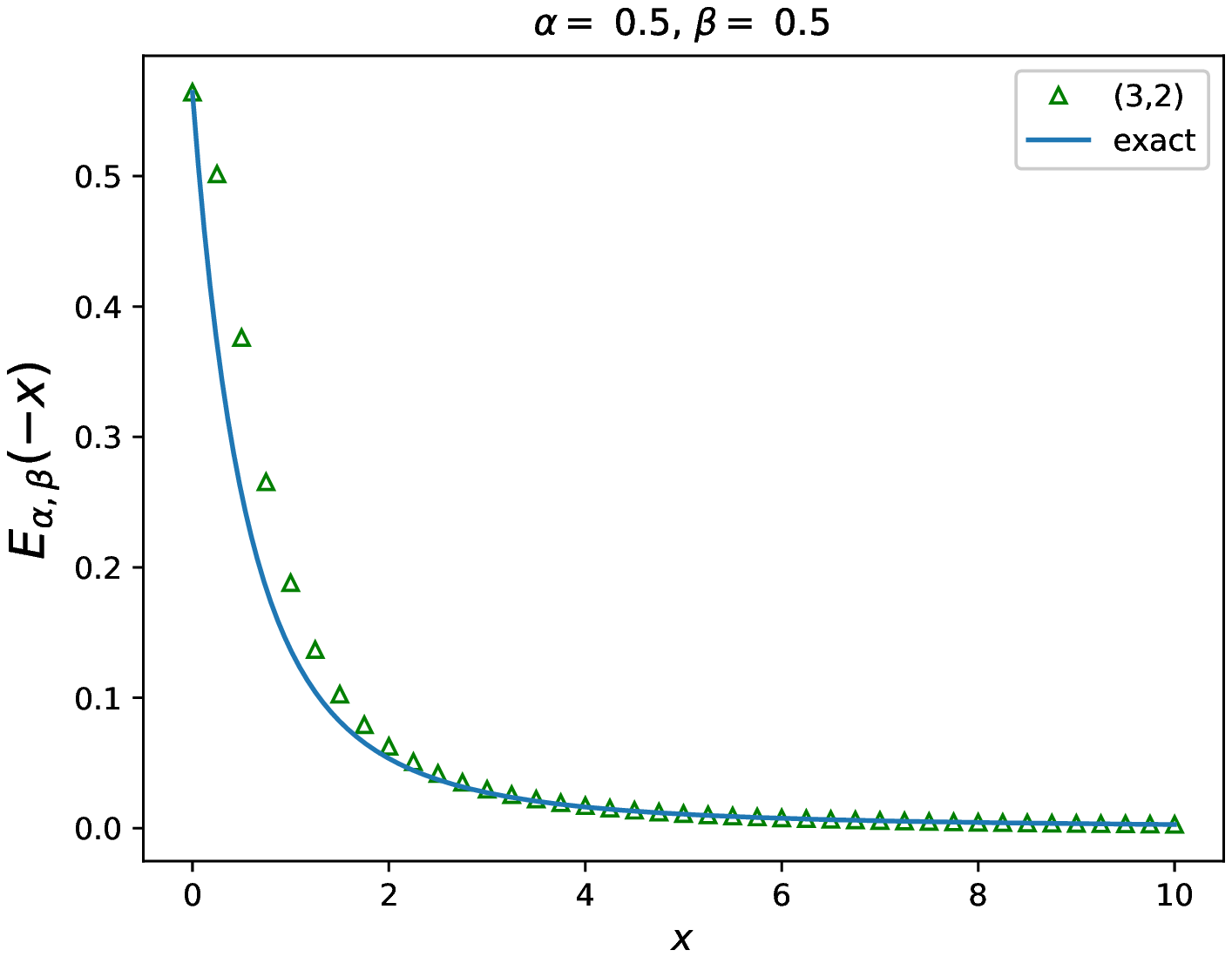}
\caption{Plots of $R^{3,2}_{\alpha,\beta}$ for $\beta = \alpha = 0.25$ (left) and $\beta = \alpha = 0.5$ (right)}
\label{fig:Zeng alpha is beta}
\end{figure}

As shown in Figures \ref{fig:Zeng small alpha beta 1}
and \ref{fig:Zeng alpha is beta}, the approximation $R^{3,2}_{\alpha,\beta}$ could be reasonable for small values of $\alpha$, however, it is not adequate otherwise.

\section{Fourth-order global Pad\'e approximants}
\label{sec:fourth-order GPA}
Fourth order global Pad\'e approximants ($\nu = 4$) correspond to the types $(m,n)$ with $m+n = 9$. 
They include the types $(5,4)$, $(6,3)$, and $(7,2)$. 
As discussed in subsection \ref{subsec:gpa construction}, the approximation 
$R_{\alpha, \beta}^{m,n}$ for $\nu = 4$ takes the form
\begin{equation}
R_{\alpha, \beta}^{m,n} (x) = 
\begin{cases}
\label{eq:Mittag Leffler gpa}
\dfrac{1}{\Gamma(\beta - \alpha)} \,
\dfrac{p_1 + p_2x + p_3 x^2 + x^3}{q_0 + q_1x + q_2 x^2 + q_3 x^3 + x^4},
\quad & \beta > \alpha,
\\ \\
\dfrac{-1}{\Gamma( - \alpha)} \,
\dfrac{\hat{p}_2  + \hat{p}_3 x + x^2}{\hat{q}_0 + \hat{q}_1x + 
\hat{q}_2 x^2 + \hat{q}_3 x^3 + x^4}, & \beta = \alpha.
\end{cases}
\end{equation}
The unknown coefficients are obtained by applying 
\eqref{eq:gpa matching} and \eqref{eq:gpa matching2}.
Below we present the systems for these coefficients for the different types.
\subsection{Coefficients of $R_{\alpha, \beta}^{5, 4}$}

For $\beta > \alpha$, the coefficients satisfy the system

\begin{equation}
\label{sys:54ab}
\begin{bmatrix}
  1 & 0  & 0 & - \dfrac{\Gamma(\beta - \alpha)}{\Gamma{(\beta)}}         
  & 0                                               & 0 & 0              \\
  0 & 1  & 0 & \dfrac{\Gamma(\beta - \alpha)}{\Gamma{(\beta + \alpha)}}  
  & -\dfrac{\Gamma(\beta - \alpha)}{\Gamma{(\beta)}}  & 0 & 0               \\
  0 & 0  & 1 & -\dfrac{\Gamma(\beta - \alpha)}{\Gamma{(\beta + 2\alpha)}} 
  & \dfrac{\Gamma(\beta - \alpha)}{\Gamma{(\beta + \alpha)}} & -\dfrac{\Gamma(\beta - \alpha)}{\Gamma{(\beta)}} & 0 \\
  0 & 0  & 0 & \dfrac{\Gamma(\beta - \alpha)}{\Gamma{(\beta + 3\alpha)}} 
  & - \dfrac{\Gamma(\beta - \alpha)}{\Gamma{(\beta + 2\alpha)}} & \dfrac{\Gamma(\beta - \alpha)}{\Gamma{(\beta + \alpha)}} 
  & - \dfrac{\Gamma{(\beta - \alpha)}}{\Gamma{(\beta)}} \\
  1 & 0  & 0 & 0 &-1 & \dfrac{\Gamma(\beta - \alpha)}{\Gamma(\beta - 2\alpha)} & 
  - \dfrac{\Gamma(\beta - \alpha)}{\Gamma(\beta - 3\alpha)}  \\
  0 & 1  & 0 & 0 & 0 &  -1  
  & \dfrac{\Gamma(\beta - \alpha)}{\Gamma(\beta - 2\alpha)}  \\
  0 & 0  & 1 & 0 & 0 & 0 & -1
\end{bmatrix}
\begin{pmatrix}
  p_1                 \\
  p_2                \\
  p_3              	\\
  q_0				  \\
  q_1 \\
  q_2 \\
  q_3
\end{pmatrix} = 
\begin{pmatrix}
  0                 \\
  0                \\
  0              	\\
  -1				  \\
  -\dfrac{\Gamma(\beta - \alpha)}{\Gamma(\beta - 4\alpha)} \\
  \dfrac{\Gamma(\beta - \alpha)}{\Gamma(\beta - 3\alpha)} \\
  -\dfrac{\Gamma(\beta - \alpha)}{\Gamma(\beta - 2\alpha)}
\end{pmatrix}.
\end{equation}

For $\beta = \alpha$ the coefficients satisfy the system
\begin{equation}
\label{sys:54aa}
\begin{bmatrix}
1  & 0 & \dfrac{\Gamma( - \alpha)}{\Gamma{(\alpha)}}         
& 0                                               & 0 & 0              \\
0  & 1 & -\dfrac{\Gamma( - \alpha)}{\Gamma{(2 \alpha)}}  
& \dfrac{\Gamma( - \alpha)}{\Gamma{(\alpha)}}  & 0 & 0               \\
0  & 0 & \dfrac{\Gamma( - \alpha)}{\Gamma{(3\alpha)}} 
& -\dfrac{\Gamma( - \alpha)}{\Gamma{(2\alpha)}} & \dfrac{\Gamma( - \alpha)}{\Gamma{(\alpha)}} & 0 \\
0  & 0 & 0 
&  -1 & -\dfrac{\Gamma( - \alpha)}{\Gamma{(-2 \alpha)}} 
& 0 \\
1  & 0 & 0 & 0 
& -1 &  \dfrac{\Gamma( - \alpha)}{\Gamma(-2\alpha)}  \\
0  & 1 & 0 & 0 &  0  
& -1 
\end{bmatrix}
\begin{pmatrix}
\hat{p}_2 \\
\hat{p}_3 \\
\hat{q}_0 \\
\hat{q}_1 \\
\hat{q}_2 \\
\hat{q}_3
\end{pmatrix} = 
\begin{pmatrix}
0                \\
0              	\\
-1				  \\
0                  \\
0                  \\
-\dfrac{\Gamma( - \alpha)}{\Gamma( - 2\alpha)}
\end{pmatrix}.
\end{equation}

\subsection{Coefficients of $R_{\alpha, \beta}^{6,3}$}	
For $\beta > \alpha$, the coefficients satisfy the system
\begin{gather}\label{sys:63ab}
\begin{bmatrix}
1 & 0  & 0 & - \dfrac{\Gamma(\beta - \alpha)}{\Gamma{(\beta)}}         
& 0                                               & 0 & 0              \\
0 & 1  & 0 & \dfrac{\Gamma(\beta - \alpha)}{\Gamma{(\beta + \alpha)}}  
& -\dfrac{\Gamma(\beta - \alpha)}{\Gamma{(\beta)}}  & 0 & 0               \\
0 & 0  & 1 & -\dfrac{\Gamma(\beta - \alpha)}{\Gamma{(\beta + 2\alpha)}} 
& \dfrac{\Gamma(\beta - \alpha)}{\Gamma{(\beta + \alpha)}} & -\dfrac{\Gamma(\beta - \alpha)}{\Gamma{(\beta)}} & 0 \\
0 & 0  & 0 & \dfrac{\Gamma(\beta - \alpha)}{\Gamma{(\beta + 3\alpha)}} 
& - \dfrac{\Gamma(\beta - \alpha)}{\Gamma{(\beta + 2\alpha)}} & \dfrac{\Gamma(\beta - \alpha)}{\Gamma{(\beta + \alpha)}} 
& - \dfrac{\Gamma{(\beta - \alpha)}}{\Gamma{(\beta)}} \\
0 & 0  & 0 & -\dfrac{\Gamma(\beta - \alpha)}{\Gamma{(\beta + 4\alpha)}} & \dfrac{\Gamma(\beta - \alpha)}{\Gamma{(\beta + 3\alpha)}} 
& -\dfrac{\Gamma(\beta - \alpha)}{\Gamma(\beta + 2\alpha)} &  \dfrac{\Gamma(\beta - \alpha)}{\Gamma(\beta + \alpha)}  \\
0 & 1  & 0 & 0 & 0 &  -1  
& \dfrac{\Gamma(\beta - \alpha)}{\Gamma(\beta - 2\alpha)}  \\
0 & 0  & 1 & 0 & 0 & 0 & -1
\end{bmatrix}
\begin{pmatrix}
p_1                 \\
p_2                \\
p_3              	\\
q_0				  \\
q_1 \\
q_2 \\
q_3
\end{pmatrix} = 
\begin{pmatrix}
0                 \\
0                \\
0              	\\
-1				  \\
\dfrac{\Gamma(\beta - \alpha)}{\Gamma(\beta)} \\
\dfrac{\Gamma(\beta - \alpha)}{\Gamma(\beta - 3\alpha)} \\
-\dfrac{\Gamma(\beta - \alpha)}{\Gamma(\beta - 2\alpha)}
\end{pmatrix},
\end{gather}
For $\beta = \alpha$, the coefficients satisfy the system
\begin{gather}\label{sys:63aa}
\begin{bmatrix}
		1  & 0 & \dfrac{\Gamma( - \alpha)}{\Gamma{(\alpha)}}         
		& 0                                               & 0 & 0              \\
		0  & 1 & -\dfrac{\Gamma( - \alpha)}{\Gamma{(2 \alpha)}}  
		& \dfrac{\Gamma( - \alpha)}{\Gamma{(\alpha)}}  & 0 & 0               \\
		0  & 0 & \dfrac{\Gamma( - \alpha)}{\Gamma{(3\alpha)}} 
		& -\dfrac{\Gamma( - \alpha)}{\Gamma{(2\alpha)}} & \dfrac{\Gamma( - \alpha)}{\Gamma{(\alpha)}} & 0 \\
		0  & 0 & -\dfrac{\Gamma( - \alpha)}{\Gamma{(4\alpha)}} 
		&  \dfrac{\Gamma( - \alpha)}{\Gamma{(3\alpha)}} & -\dfrac{\Gamma( - \alpha)}{\Gamma{(2 \alpha)}} 
		& \dfrac{\Gamma{( - \alpha)}}{\Gamma{(\alpha)}} \\
		1  & 0 & 0 & 0 
		& -1 &  \dfrac{\Gamma( - \alpha)}{\Gamma(-2\alpha)}  \\
		0  & 1 & 0 & 0 &  0  
		& -1 
\end{bmatrix}
\begin{pmatrix}
	\hat{p}_2 \\
	\hat{p}_3 \\
	\hat{q}_0 \\
	\hat{q}_1 \\
	\hat{q}_2 \\
	\hat{q}_3
\end{pmatrix} = 
\begin{pmatrix}
		0                \\
		0              	\\
		-1				  \\
		0                  \\
		\dfrac{\Gamma( - \alpha)}{\Gamma( - 3\alpha)} \\
		-\dfrac{\Gamma( - \alpha)}{\Gamma( - 2\alpha)}
\end{pmatrix}.
\end{gather}

\subsection{Coefficients of $R_{\alpha, \beta}^{7, 2}$}
For $\beta > \alpha$, the coefficients satisfy the system
\begin{gather}\label{sys:72ab}
\begin{bmatrix}
1 & 0  & 0 & - \dfrac{\Gamma(\beta - \alpha)}{\Gamma{(\beta)}}         
& 0                                               & 0 & 0              \\
0 & 1  & 0 & \dfrac{\Gamma(\beta - \alpha)}{\Gamma{(\beta + \alpha)}}  
& -\dfrac{\Gamma(\beta - \alpha)}{\Gamma{(\beta)}}  & 0 & 0               \\
0 & 0  & 1 & -\dfrac{\Gamma(\beta - \alpha)}{\Gamma{(\beta + 2\alpha)}} 
& \dfrac{\Gamma(\beta - \alpha)}{\Gamma{(\beta + \alpha)}} & -\dfrac{\Gamma(\beta - \alpha)}{\Gamma{(\beta)}} & 0 \\
0 & 0  & 0 & \dfrac{\Gamma(\beta - \alpha)}{\Gamma{(\beta + 3\alpha)}} 
& - \dfrac{\Gamma(\beta - \alpha)}{\Gamma{(\beta + 2\alpha)}} & \dfrac{\Gamma(\beta - \alpha)}{\Gamma{(\beta + \alpha)}} 
& - \dfrac{\Gamma{(\beta - \alpha)}}{\Gamma{(\beta)}} \\
0 & 0  & 0 & -\dfrac{\Gamma(\beta - \alpha)}{\Gamma{(\beta + 4\alpha)}} & \dfrac{\Gamma(\beta - \alpha)}{\Gamma{(\beta + 3\alpha)}} 
& -\dfrac{\Gamma(\beta - \alpha)}{\Gamma(\beta + 2\alpha)} &  \dfrac{\Gamma(\beta - \alpha)}{\Gamma(\beta + \alpha)}  \\
0 & 0  & 0 & \dfrac{\Gamma(\beta - \alpha)}{\Gamma{(\beta + 5\alpha)}} & -\dfrac{\Gamma(\beta - \alpha)}{\Gamma{(\beta + 4\alpha)}} 
& \dfrac{\Gamma(\beta - \alpha)}{\Gamma(\beta + 3\alpha)} &  -\dfrac{\Gamma(\beta - \alpha)}{\Gamma(\beta + 2\alpha)}  \\
0 & 0  & 1 & 0 & 0 & 0 & -1
\end{bmatrix}
\begin{pmatrix}
p_1  \\
p_2   \\
p_3 \\
q_0	\\
q_1 \\
q_2 \\
q_3
\end{pmatrix} = 
\begin{pmatrix}
0  \\
0   \\
0   \\
-1	\\
\dfrac{\Gamma(\beta - \alpha)}{\Gamma(\beta)} \\
-\dfrac{\Gamma(\beta - \alpha)}{\Gamma(\beta + \alpha)} \\
-\dfrac{\Gamma(\beta - \alpha)}{\Gamma(\beta - 2\alpha)}
\end{pmatrix}.
\end{gather}

For $\beta = \alpha$, the coefficients satisfy the system
\begin{gather}\label{sys:72aa}
\begin{bmatrix}
1  & 0 & \dfrac{\Gamma( - \alpha)}{\Gamma{(\alpha)}}         
& 0                                               & 0 & 0              \\
0  & 1 & -\dfrac{\Gamma( - \alpha)}{\Gamma{(2 \alpha)}}  
& \dfrac{\Gamma( - \alpha)}{\Gamma{(\alpha)}}  & 0 & 0               \\
0  & 0 & \dfrac{\Gamma( - \alpha)}{\Gamma{(3\alpha)}} 
& -\dfrac{\Gamma( - \alpha)}{\Gamma{(2\alpha)}} & \dfrac{\Gamma( - \alpha)}{\Gamma{(\alpha)}} & 0 \\
0  & 0 & -\dfrac{\Gamma( - \alpha)}{\Gamma{(4\alpha)}} 
&  \dfrac{\Gamma( - \alpha)}{\Gamma{(3\alpha)}} & -\dfrac{\Gamma( - \alpha)}{\Gamma{(2 \alpha)}} 
& - \dfrac{\Gamma{( - \alpha)}}{\Gamma{(\alpha)}} \\
0  & 0 & \dfrac{\Gamma( - \alpha)}{\Gamma{(5\alpha)}} 
&  -\dfrac{\Gamma( - \alpha)}{\Gamma{(4\alpha)}} & \dfrac{\Gamma( - \alpha)}{\Gamma{(3 \alpha)}} 
& - \dfrac{\Gamma{( - \alpha)}}{\Gamma{(2\alpha)}} \\
0  & 1 & 0 & 0 &  0  
& -1 
\end{bmatrix}
\begin{pmatrix}
	\hat{p}_2 \\
	\hat{p}_3 \\
	\hat{q}_0 \\
	\hat{q}_1 \\
	\hat{q}_2 \\
	\hat{q}_3
\end{pmatrix} = 
\begin{pmatrix}
0  \\
0  	\\
-1  \\
0  \\
\dfrac{\Gamma( - \alpha)}{\Gamma( - 3\alpha)} \\
-\dfrac{\Gamma( - \alpha)}{\Gamma( - 2\alpha)}
\end{pmatrix}.
\end{gather}

\begin{remark}
Although the type $(m,n)$ of the fourth order global Pad\'e approximant of $E_\alpha(-x)$ by Ingo \etal \cite{Ingo:2017} is not given, 
based on their construction, we expect the type to be a special case of one of the approximants above when $\beta = 1$.
\end{remark}
\section{Partial fraction decomposition}
\label{sec:partial-frac}

Partial fraction decomposition provides an efficient form for evaluating rational functions. 
In a recent work by \cite{Bertaccini:2019}, the efficiency of using partial fraction decomposition for computing functions of matrices is discussed.
This efficiency is indisputable when the poles are complex conjugates and the argument is a matrix.

Unlike the the Pad\'e approximations for the exponential function, 
the poles of $R^{m,n}_{\alpha,\beta}$ are functions of $\alpha$ and $\beta$.
Fortunately, through direct calculations, one can show that for most $(\alpha,\beta) \in \mathcal{A}$, the poles of 
$R^{m,n}_{\alpha,\beta}$ are complex conjugates. Next, we explore the partial fraction decomposition of these approximations.

\subsection{Decomposition of second-order global Pad\'e approximant}

The second-order global Pad\'e approximant $R^{3,2}_{\alpha, \beta}$ admits the partial fraction decomposition
$$
R^{3,2}_{\alpha, \beta}(x) = \frac{c_1}{x - r_1} + \frac{c_2}{x - r_2},
$$ 
where 
$$
r_1 = \frac{-q_1 + \sqrt{q_1^2 - 4q_0}}{2}\text{, }\qquad  
r_2 = \frac{-q_1 - \sqrt{q_1^2 - 4q_0}}{2},
$$ and
$$
c_1 = \frac{p_1 - r_1}{r_2 - r_1}\text{, } \qquad  
c_2 = \frac{p_1 - r_2}{r_1 - r_2}.
$$
We can verify numerically that for $(\alpha,\beta) \in \mathcal{A}$ we have $q_1^2 - 4 q_0 < 0$ 
which imply that $r_2 = \bar{r}_1$ and $c_2 = \bar{c}_1$.
As a result, we can write
\begin{equation}\label{eq:partial-frac second order}
R^{3,2}_{\alpha, \beta}(x) = 2\Re \left[\frac{c_1}{x - r_1}\right].
\end{equation}

\bigskip

\subsection{Decomposition of the fourth-order global Pad\'e approximants}

The partial fraction decomposition for $R_{\alpha, \beta}^{m, n}$, $(m,n) = (5,4), (6,3), (7,2)$, takes the form
\begin{equation}
\label{eq:partial-frac fourth order}
R_{\alpha, \beta}^{m, n}(x) = \frac{c_1}{x - r_1} + \frac{c_2}{x - r_2} + 
\frac{c_3}{x - r_3} + \frac{c_4}{x - r_4}.
\end{equation}
Empirically, for $(\alpha,\beta) \in \mathcal{A}$, these poles are complex conjugates.
If we let $r_3 = \bar{r}_1$, $r_4 = \bar{r}_2$, $c_3 = \bar{c}_1$, and 
$c_4 = \bar{c}_2$, then the partial fraction decomposition can be written as
\begin{equation}
\label{eq:partial-frac conjugate fourth order}
R_{\alpha, \beta}^{m, n}(x) = 2 \Re \left[\frac{c_1}{x - r_1} \right] 
+ 2 \Re \left[ \frac{c_2}{x - r_2}\right].
\end{equation}
Computing of the poles and weights is outlined in the following algorithm.
\begin{algorithm}
\caption{Poles and weights for partial fraction decomposition of fourth-order $R_{\alpha, \beta}^{m,n}$}
\begin{algorithmic}
\item 
\textbf{Step 1}
\\  Specify $m \text{, } n$, $\alpha$, $\beta$.
\\  Obtain $p_i$, $q_i$ by solving the corresponding system.
\item 
\textbf{Step 2} 
\\ Use the obtained coefficients $p_i, \, q_i$ to find the weights and poles: 
\\ If $\beta > \alpha$
\\ \hspace{0.5 in} Matlab: residue$([1,p_3,p_2,p_1],[1,q_3,q_2,q_1,q_0])$,
\\ \hspace{0.5 in} Python: scipy.signal.residue$([1,p_3,p_2,p_1],[1,q_3,q_2,q_1,q_0])$.
\\ If $\alpha = \beta$
\\ \hspace{0.5 in} Matlab: residue$([1,p_3,p_2],[1,q_3,q_2,q_1,q_0])$,
\\ \hspace{0.5 in} Python: scipy.signal.residue$([1,p_3,p_2],[1,q_3,q_2,q_1,q_0])$.
\end{algorithmic}	
\end{algorithm}

\section{Inverse Mittag-Leffler Function}
\label{sec:Inv-mlf}
The invertibility of $E_{\alpha,\beta}(-x)$, $x > 0$, follows from the complete monotonicity property of $E_{\alpha,\beta}$.
As shown in \cite{Gorenflo:2014}, this function is completely monotone 
if and only if $0<\alpha\le 1$ and $\beta \ge \alpha$.
Since $E_{\alpha,\beta}(0) = 1/\Gamma(\beta)$ and $\lim_{x\to\infty}E_{\alpha,\beta}(-x) = 0$, 
then for $0 < \alpha \le 1$ and $\beta \ge \alpha$,
the inverse function $-L_{\alpha,\beta}$ of $E_{\alpha,\beta}(-x)$, $x>0$, is the function 
$$
-L_{\alpha,\beta} : (0, 1/\Gamma(\beta)] \rightarrow [0, \infty),
$$
such that
\begin{equation}
-L_{\alpha,\beta}(x) = y \text{ iff } x = E_{\alpha, \beta}(-y).
\end{equation}

The inverse function can be approximated by inverting the approximations 
$R^{m,n}_{\alpha,\beta}$ of $E_{\alpha,\beta}(-x)$,
$$
E_{\alpha, \beta}(-x) \approx R^{m,n}_{\alpha,\beta} = \dfrac{1}{s_{\alpha, \beta}(x)}\dfrac{p(x)}{q(x)}.
$$
This is equivalent to computing the positive root $r^+$ of the equation
\begin{equation}
\label{eq:inverse gpa}
s_{\alpha, \beta}(x)q(x)y - p(x) = 0,
\end{equation}
where $y = E_{\alpha, \beta}(-r^+)$. 
Zeng and Chen in \cite{Zeng2015} solved a quadratic equation of the form \eqref{eq:inverse gpa} for $R_{\alpha,\beta}^{3,2}(x)$ 
to obtain the approximation
\begin{align}
-L_{\alpha, \beta}(y) &\approx \dfrac{1}{2\Gamma(\beta - \alpha)y} - \frac{q_1}{2} \nonumber \\
& +  \sqrt{\left(\frac{q_1}{2}-\frac{1}{2\Gamma(\beta - \alpha)y}\right)^2 - q_0\left(1-\frac{1}{\Gamma(\beta)y}\right)}, \quad \beta > \alpha,
\end{align}
and
\begin{align}
-L_{\alpha, \alpha}(y) &\approx - \dfrac{\Gamma(1-\alpha)}{\Gamma(1 - 2\alpha)y} \nonumber \\
& +  \sqrt{\frac{\Gamma^2(1-\alpha)}{\Gamma^2(1-2\alpha)y^2} 
	- \frac{1+\alpha}{1-\alpha}\left(1 - \frac{1}{\Gamma(\alpha)y}\right)}.
\end{align}

Approximating the inverse using our fourth order approximants involves 
finding the positive root of the fourth degree polynomial in \eqref{eq:inverse gpa}, which takes the form
\begin{align}
c_{\alpha,\beta} \, y \, q_0 - p_1 + (c_{\alpha,\beta} y q_1 - p_2) x + (c_{\alpha,\beta} y q_2 - p_3) x^2 &
\nonumber \\ \label{eq:inverse gpa fourth order}
+ (c_{\alpha,\beta} y q_3 - 1) x^3 + c_{\alpha,\beta} y x^4 = 0, & \quad c_{\alpha,\beta} = \Gamma(\beta - \alpha)\text{ for } \beta > \alpha,
\\ \nonumber \\
c_{\alpha} y\hat{q}_0 + \hat{p}_2 + (c_\alpha y \hat{q}_1 +  \hat{p}_3) x + (c_\alpha y \hat{q}_2 + 1) x^2 &
\nonumber \\ \label{eq:inverse gpa fourth order aa}
+ c_\alpha y \hat{q}_3 x^3 + c_\alpha y x^4 = 0, & \quad c_{\alpha} = \Gamma( - \alpha)\text{ for } \beta = \alpha.
\end{align}
 
It is quite tedious to solve \eqref{eq:inverse gpa fourth order} and 
\eqref{eq:inverse gpa fourth order aa} analytically. 
Alternatively, the following algorithm can be employed. 
One can verify numerical that each of \eqref{eq:inverse gpa fourth order} and \eqref{eq:inverse gpa fourth order aa} has a unique positive root.

\begin{algorithm}
\caption{Approximation of the inverse MLF}
\begin{algorithmic}
\item \\
\textbf{Step 1}
\\  
Specify $m \text{, } n$, $\alpha$, $\beta$, $y$.

\bigskip
\item 
\textbf{Step 2}  
\\ \\
If $\beta > \alpha$
\\
Compute the coefficient $c_{\alpha,\beta} = \Gamma(\beta - \alpha)$
\\
Obtain $p_i$, $q_i$ by solving the corresponding linear system.
\\ 
Find the roots of the polynomial \eqref{eq:inverse gpa fourth order}
\\ 
Then $-L_{\alpha, \beta}(y)$ is the unique positive root.
\\ \\
If $\alpha = \beta$
\\ Compute the coefficient $c_\alpha = \Gamma( - \alpha)$
\\
Obtain $\hat{p}_i$, $\hat{q}_i$ by solving the corresponding linear system.
\\ 
Find the roots of the polynomial \eqref{eq:inverse gpa fourth order aa}
\\ 
Then $-L_{\alpha, \alpha}(y)$ is the unique positive root.
\end{algorithmic}	
\end{algorithm}

\section{Performance and comparisons of the approximants}
\label{sec:comparison}

\begin{figure}
	\centering
	\includegraphics[width=0.49\textwidth]{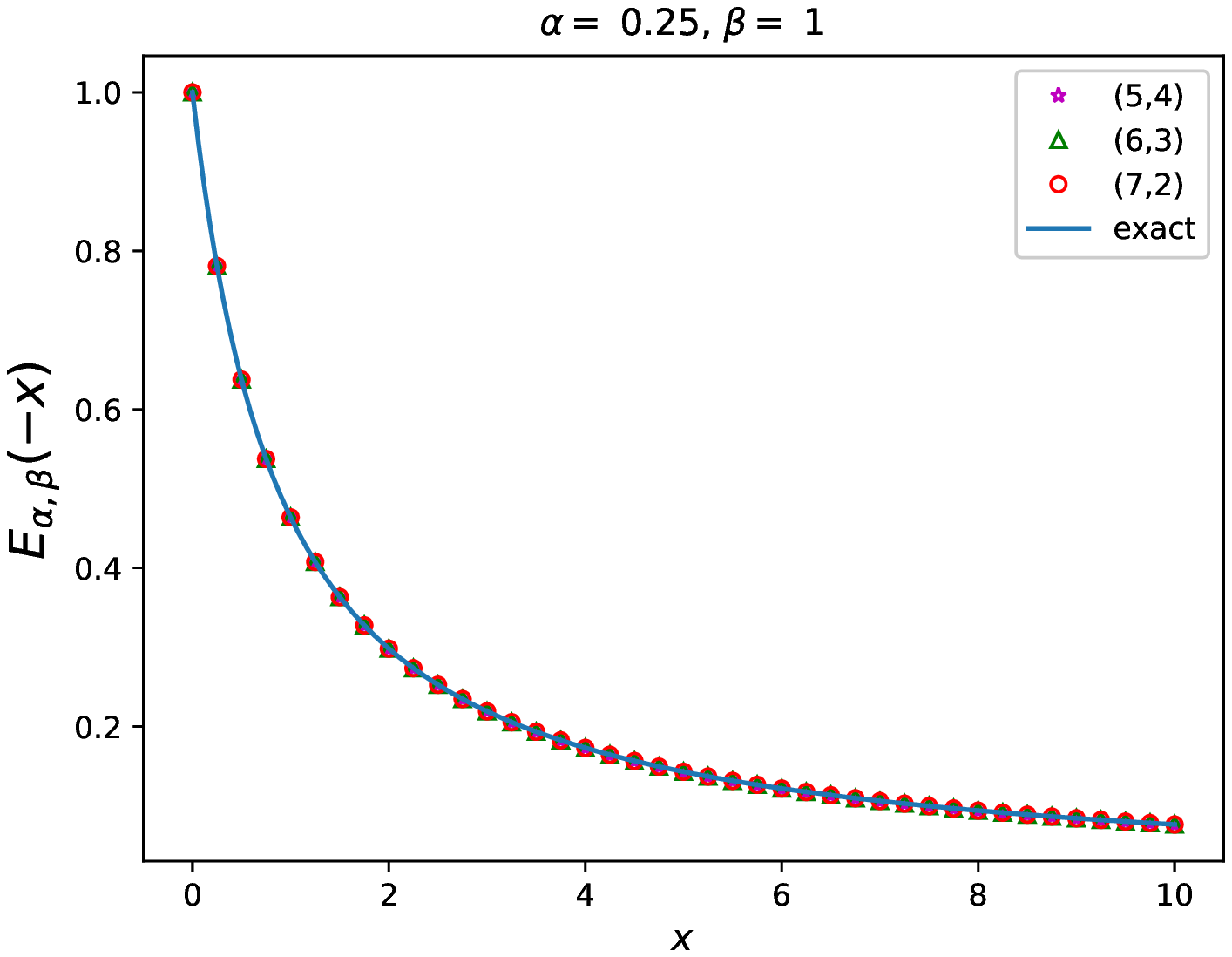}
	\includegraphics[width=0.49\textwidth]{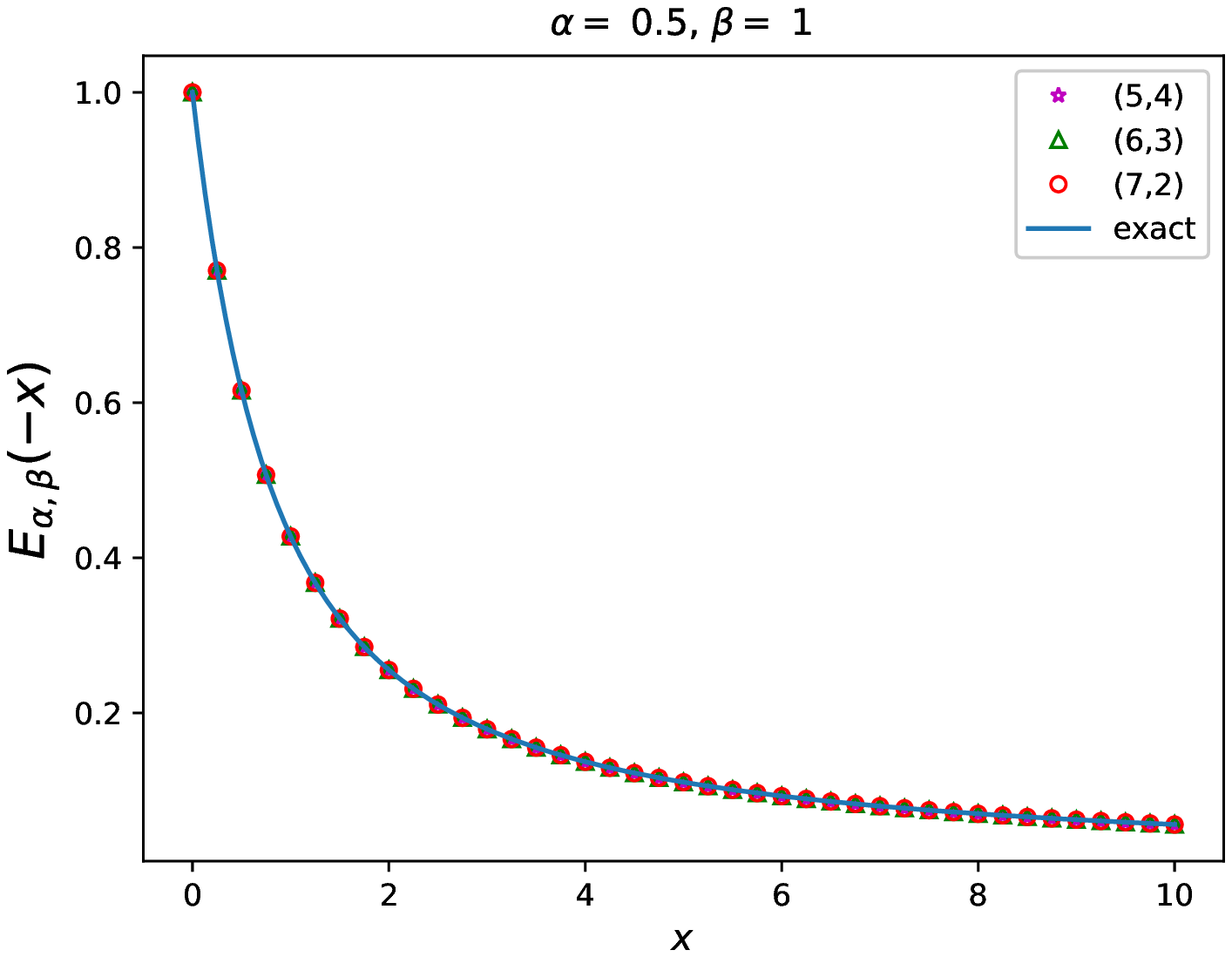} \\
	\includegraphics[width=0.49\textwidth]{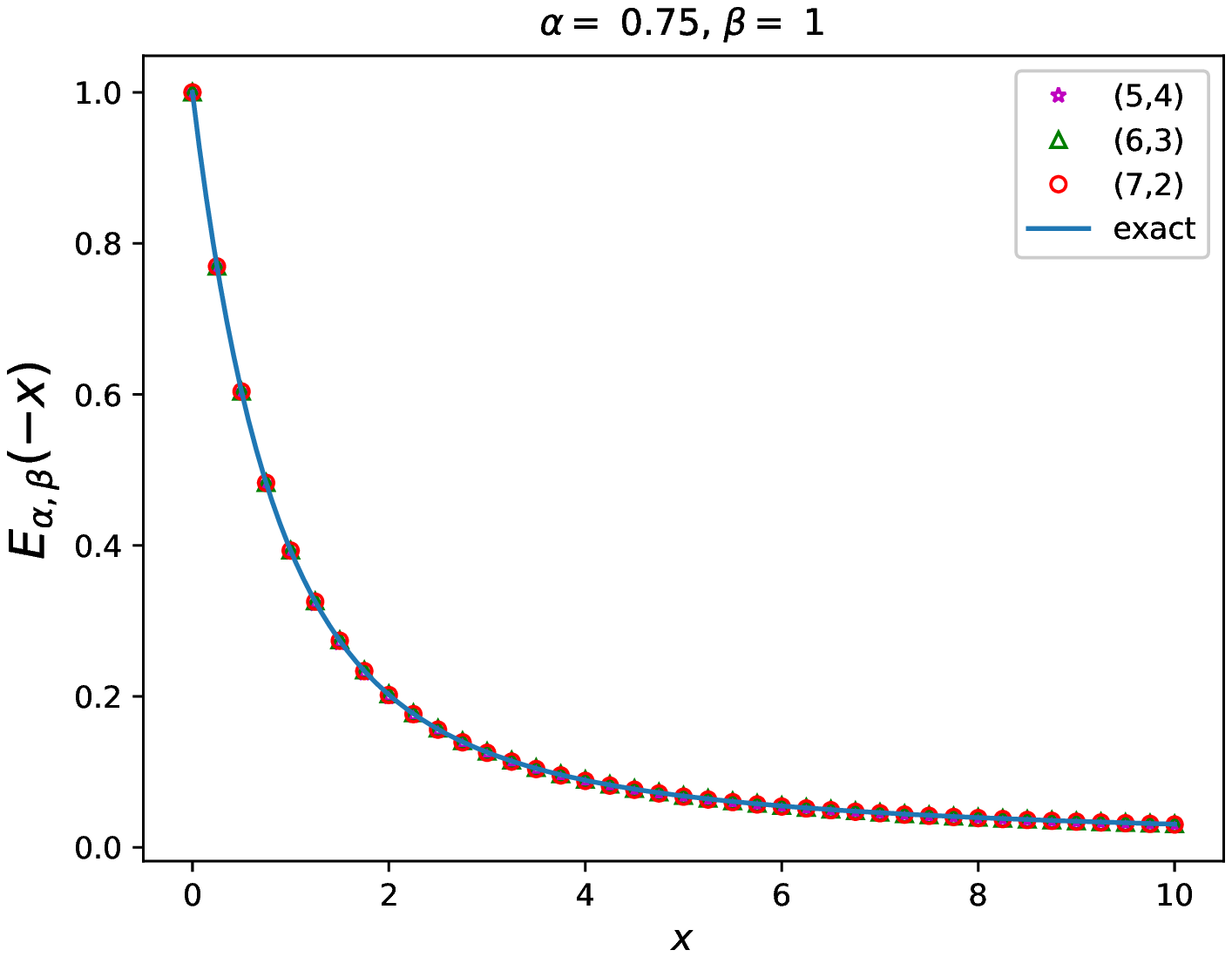}
	\includegraphics[width=0.49\textwidth]{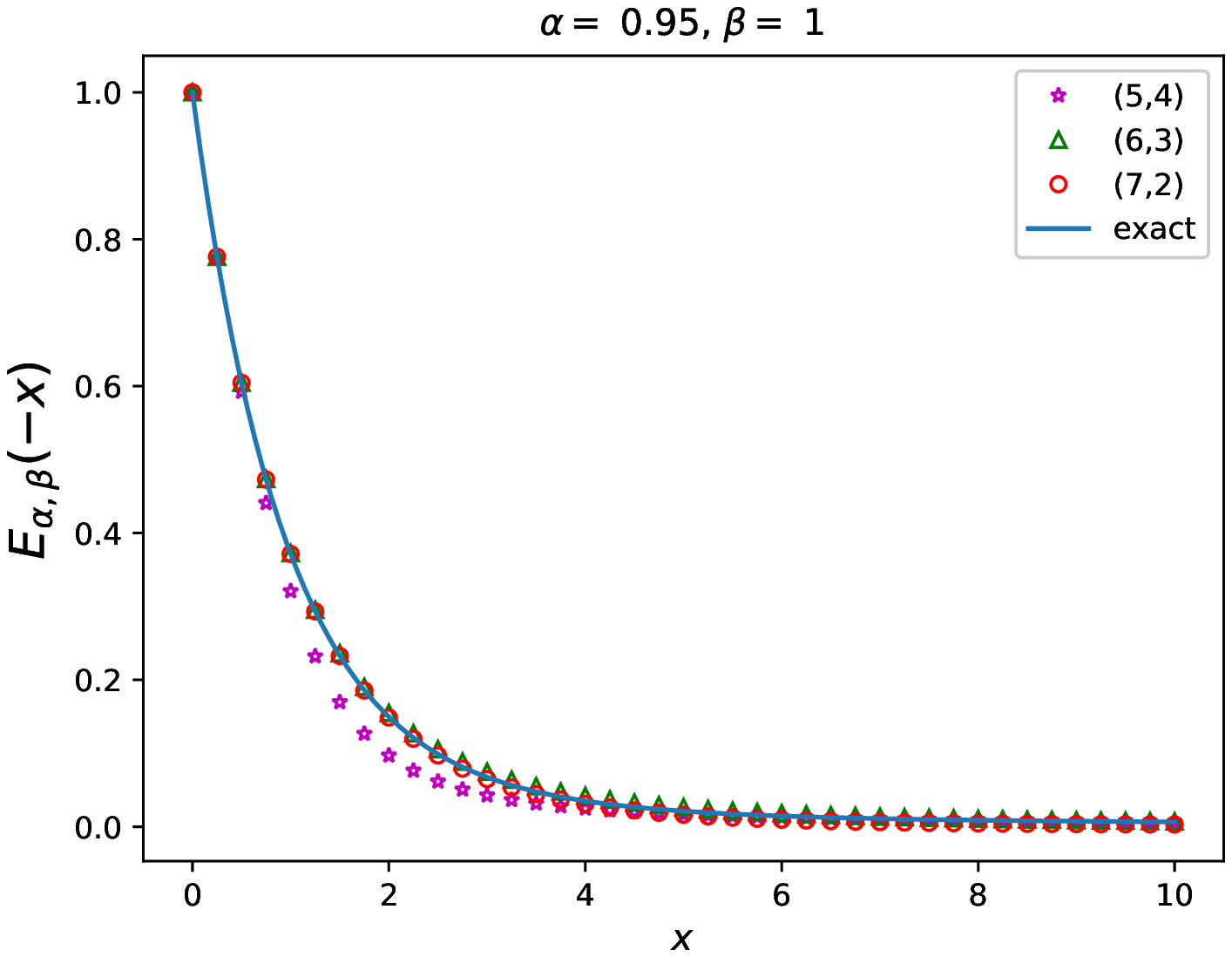} 
	\caption{Comparison of the fourth order approximants; 
		$R_{\alpha,\beta}^{5,4}$, $R_{\alpha,\beta}^{6,3}$, $R_{\alpha,\beta}^{7,2}$
		for $\beta = 1$}
	\label{fig:fourth order beta 1}
\end{figure}

\begin{figure}
	\centering
	\includegraphics[width=0.49\textwidth]{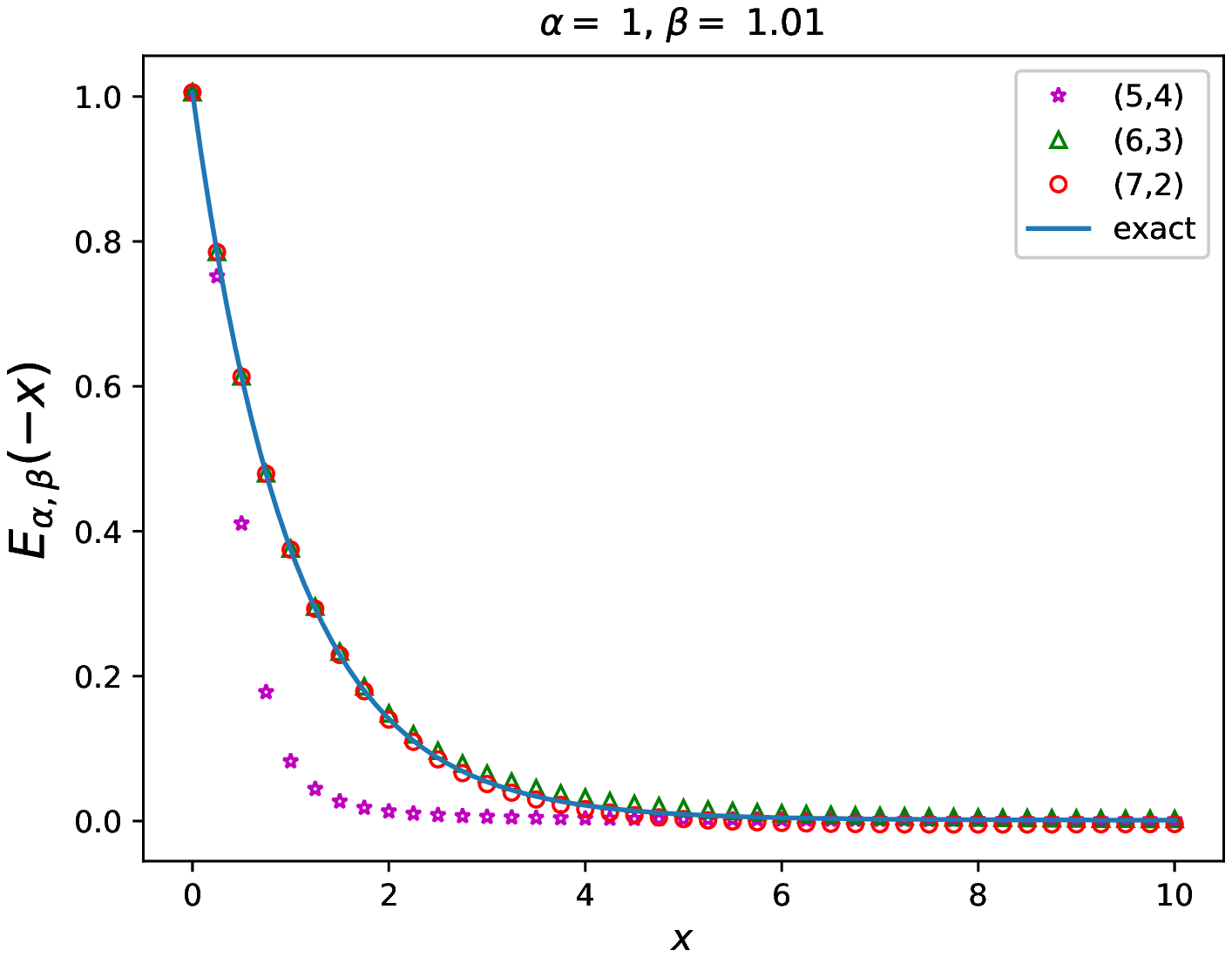}
	\includegraphics[width=0.49\textwidth]{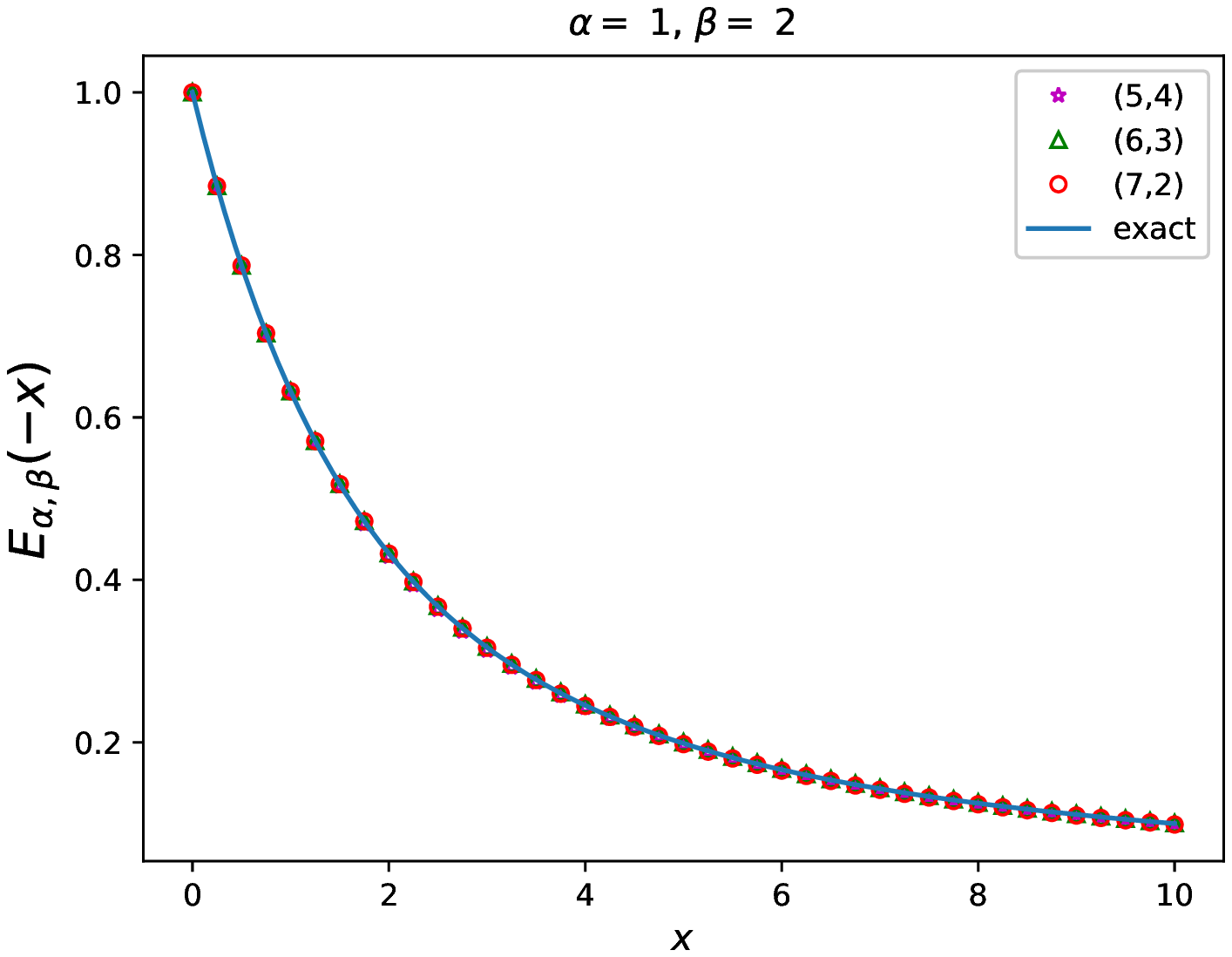}
	\caption{Comparison of the fourth order approximants;
		$R_{1,\beta}^{5,4}$, $R_{1,\beta}^{6,3}$, $R_{1,\beta}^{7,2}$
		for $\beta = 1.01$ (left) and $\beta = 2$ (right)}
	\label{fig:fourth order alpha 1}
\end{figure}

\begin{figure}
\centering
\includegraphics[width=0.49\textwidth]{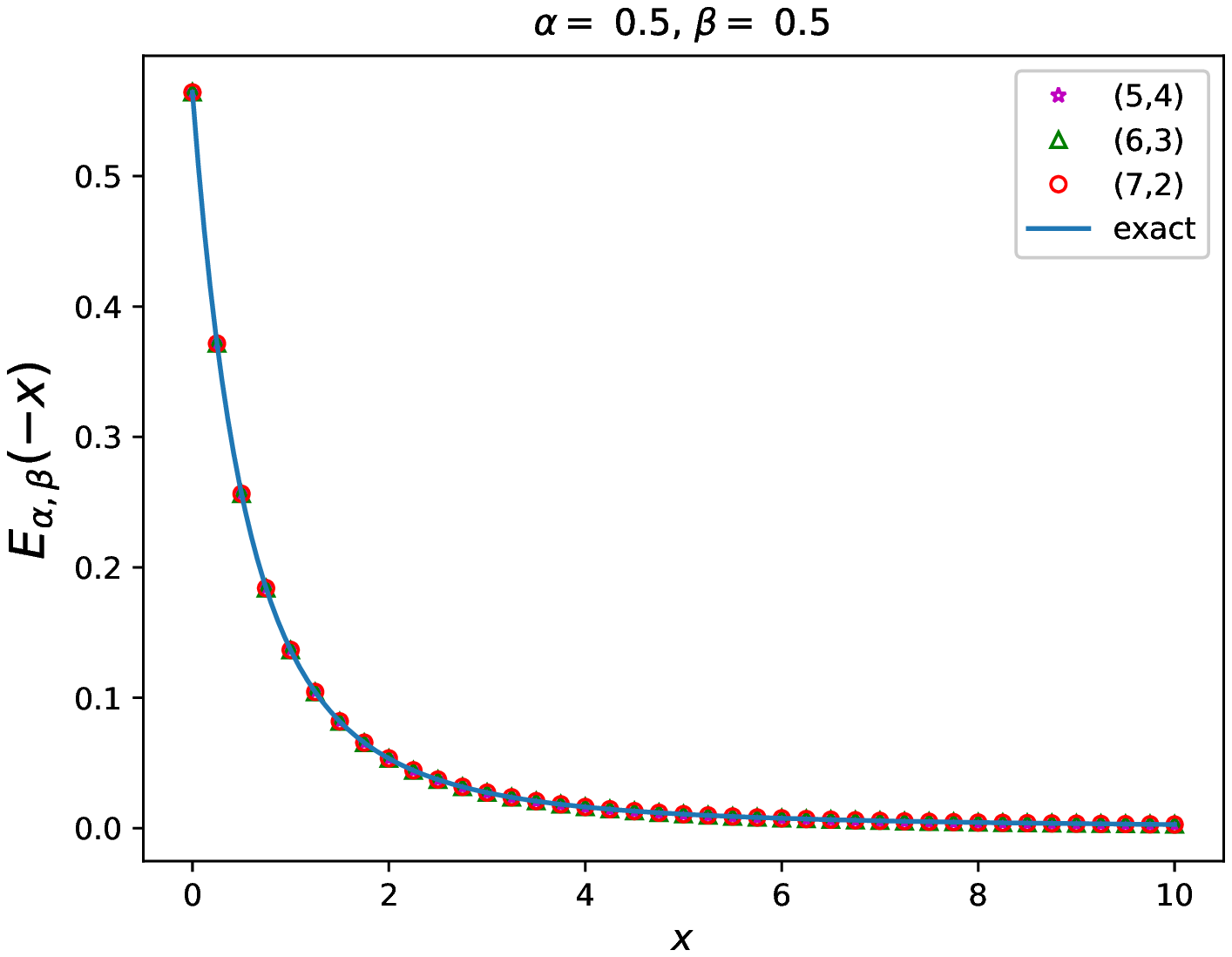}
\includegraphics[width=0.49\textwidth]{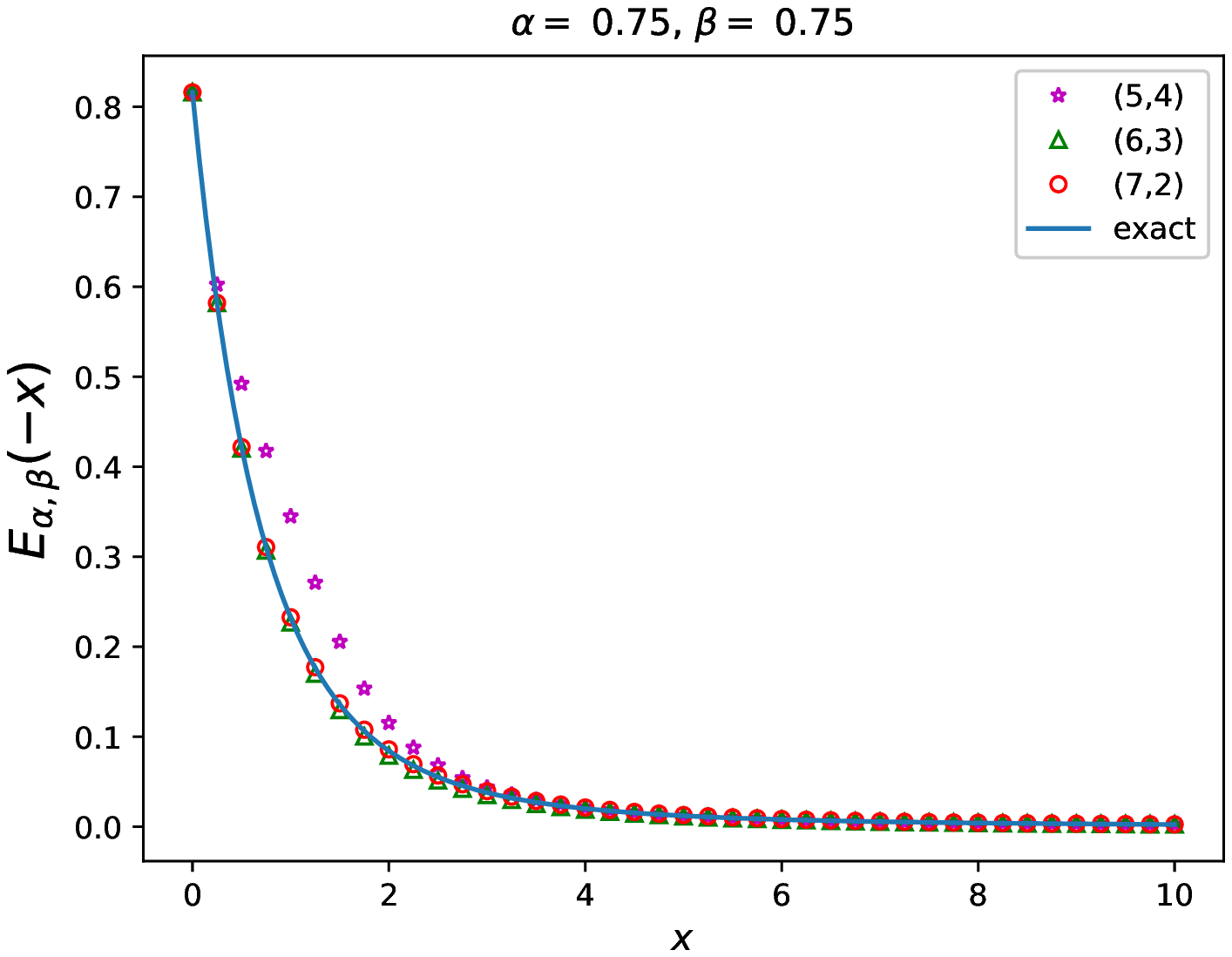}
\caption{Comparison of the fourth order approximants; 
	$R_{\alpha,\alpha}^{5,4}$, $R_{\alpha,\alpha}^{6,3}$, $R_{\alpha,\alpha}^{7,2}$
	for $\alpha = 0.5$ (left) and $\alpha = 0.75$ (right)}
\label{fig:fourth order alpha is beta}
\end{figure}

\begin{figure}
	\centering
	\includegraphics[width=0.49\textwidth]{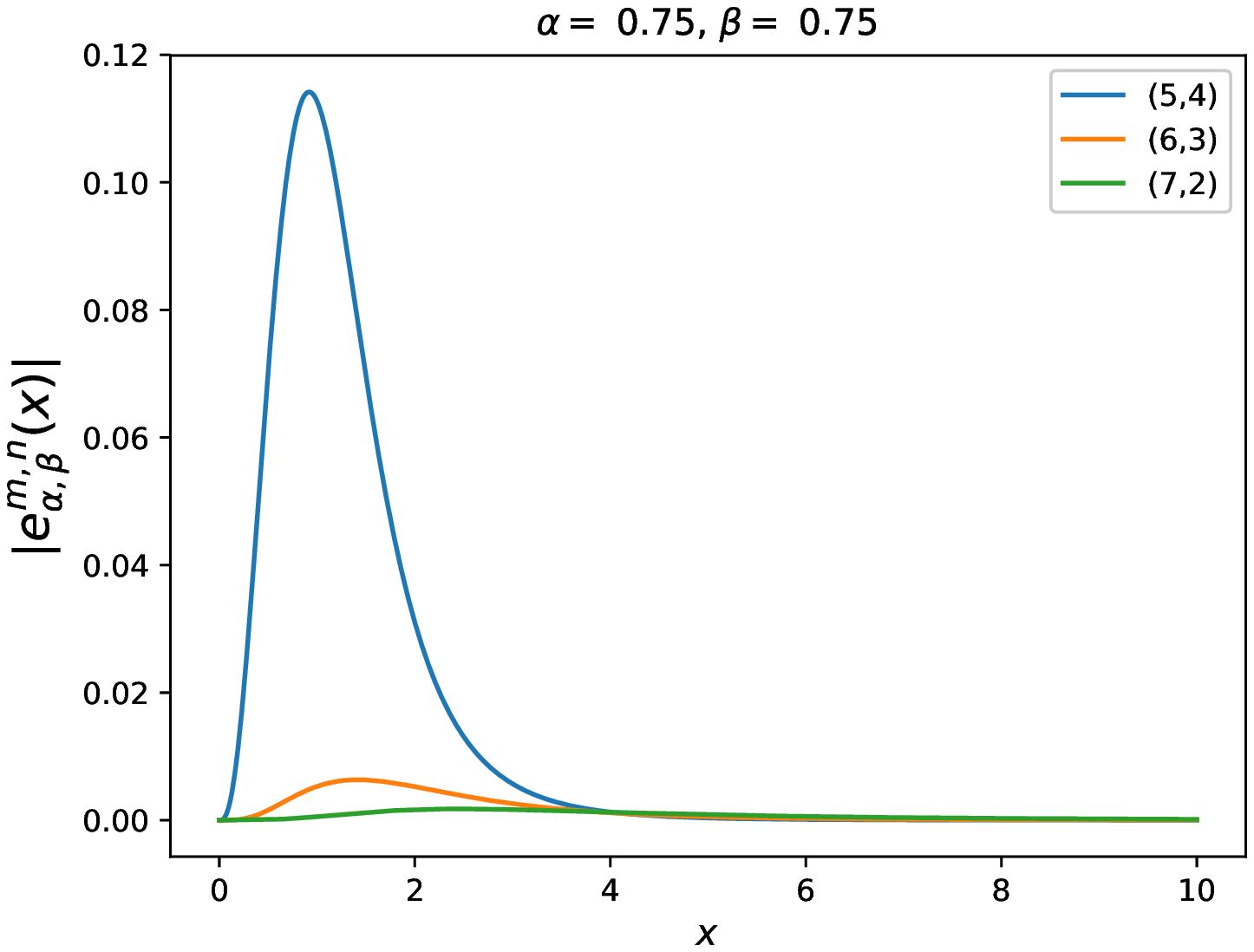}
	\includegraphics[width=0.49\textwidth]{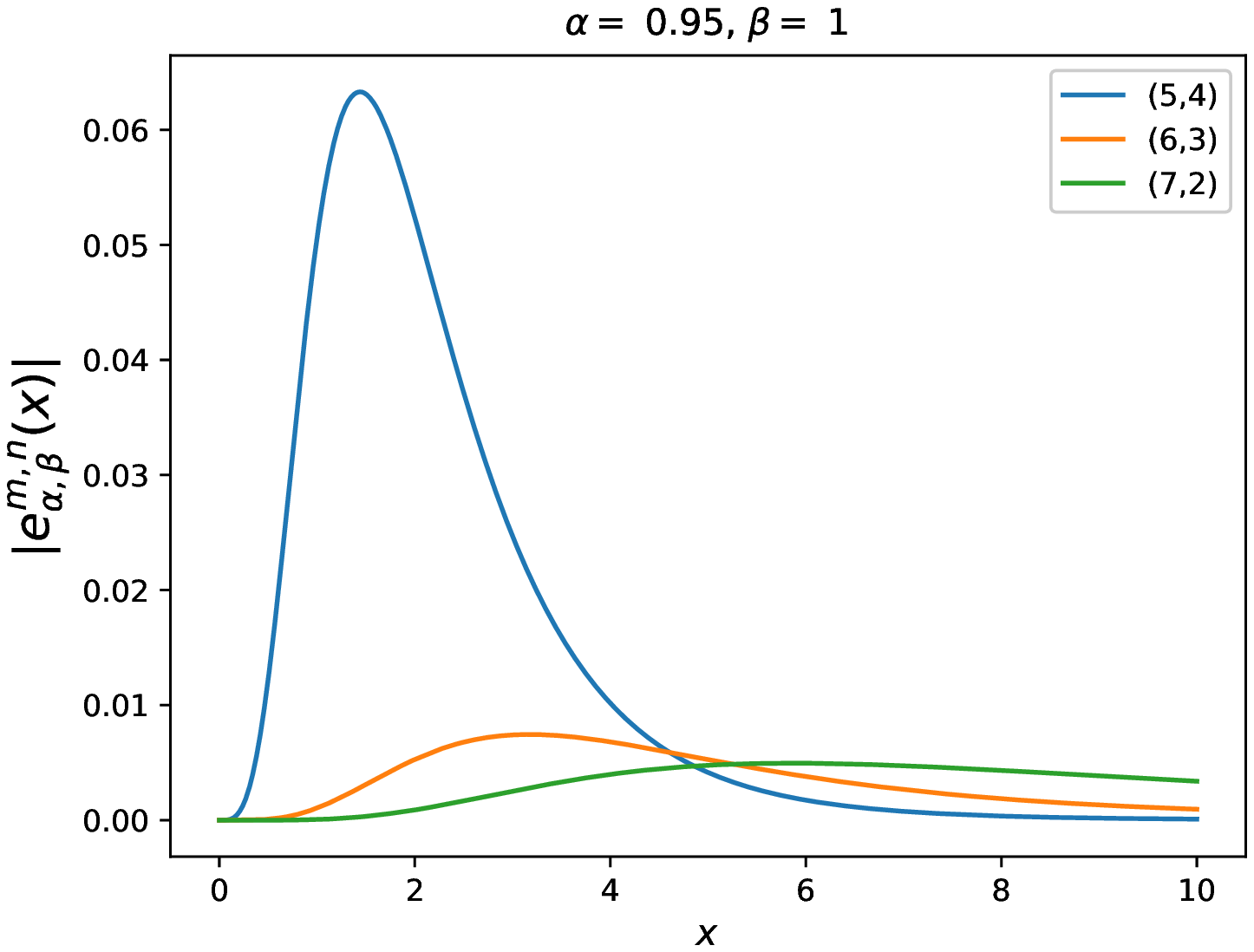} 
	\\
	\includegraphics[width=0.49\textwidth]{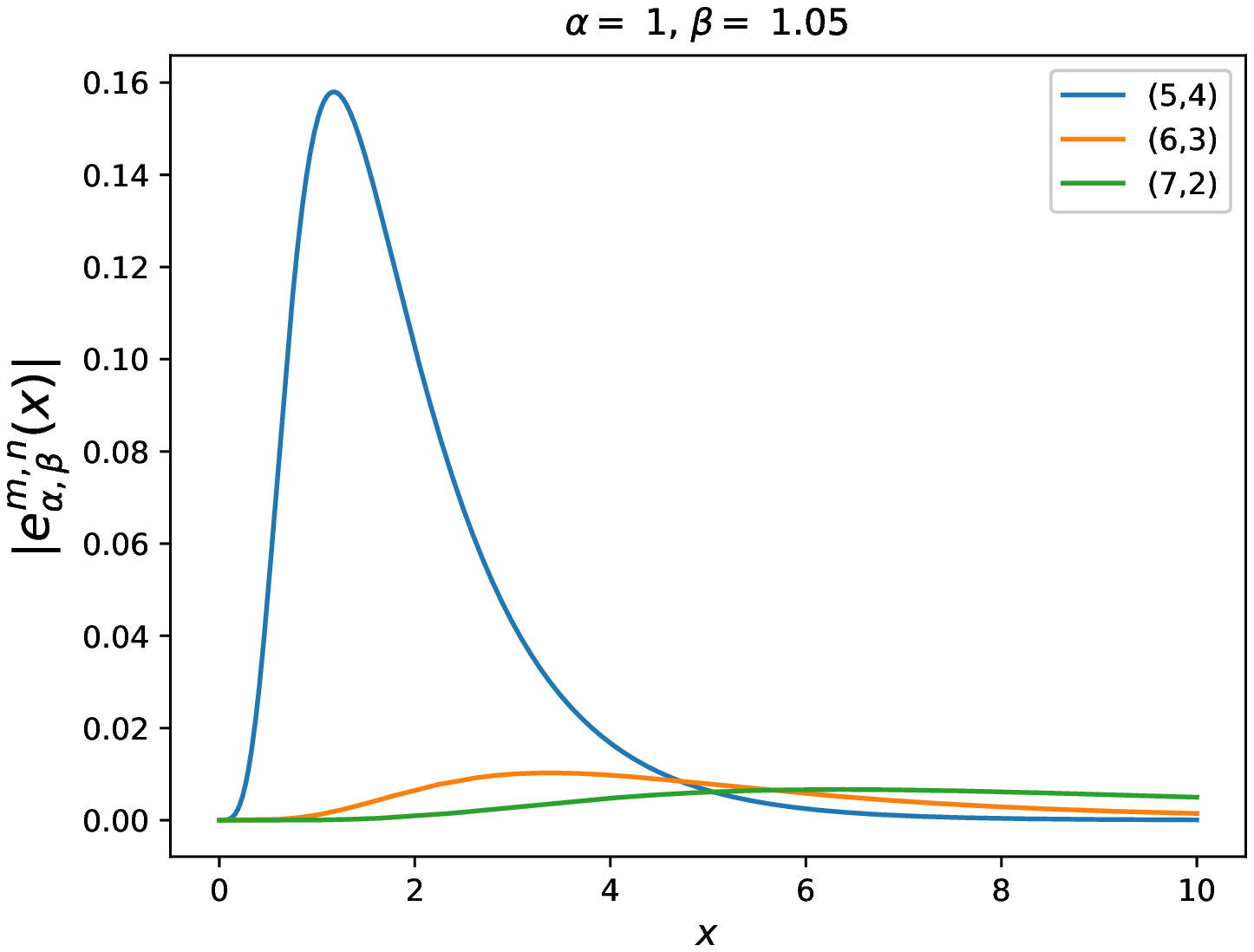}
	\includegraphics[width=0.49\textwidth]{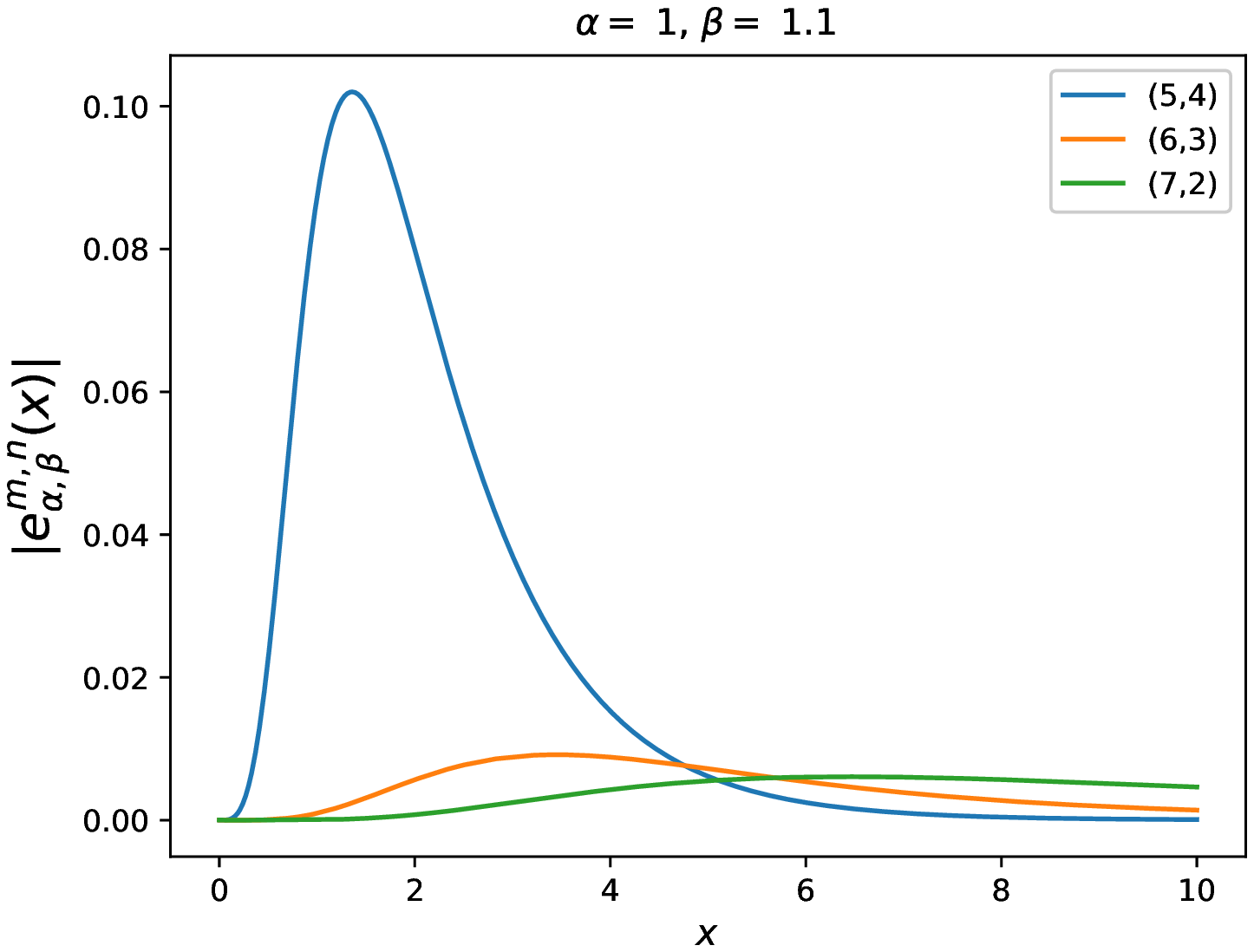}
	\caption{Plots of approximation errors $\abs{e^{m,n}_{\alpha, \beta}(x)}$ 
		for the fourth order approximants (7,2), (6,3), (5,4)}
	\label{fig:fourth order approximation errors}
\end{figure}

\begin{figure}
\centering
\includegraphics[width=0.49\textwidth]{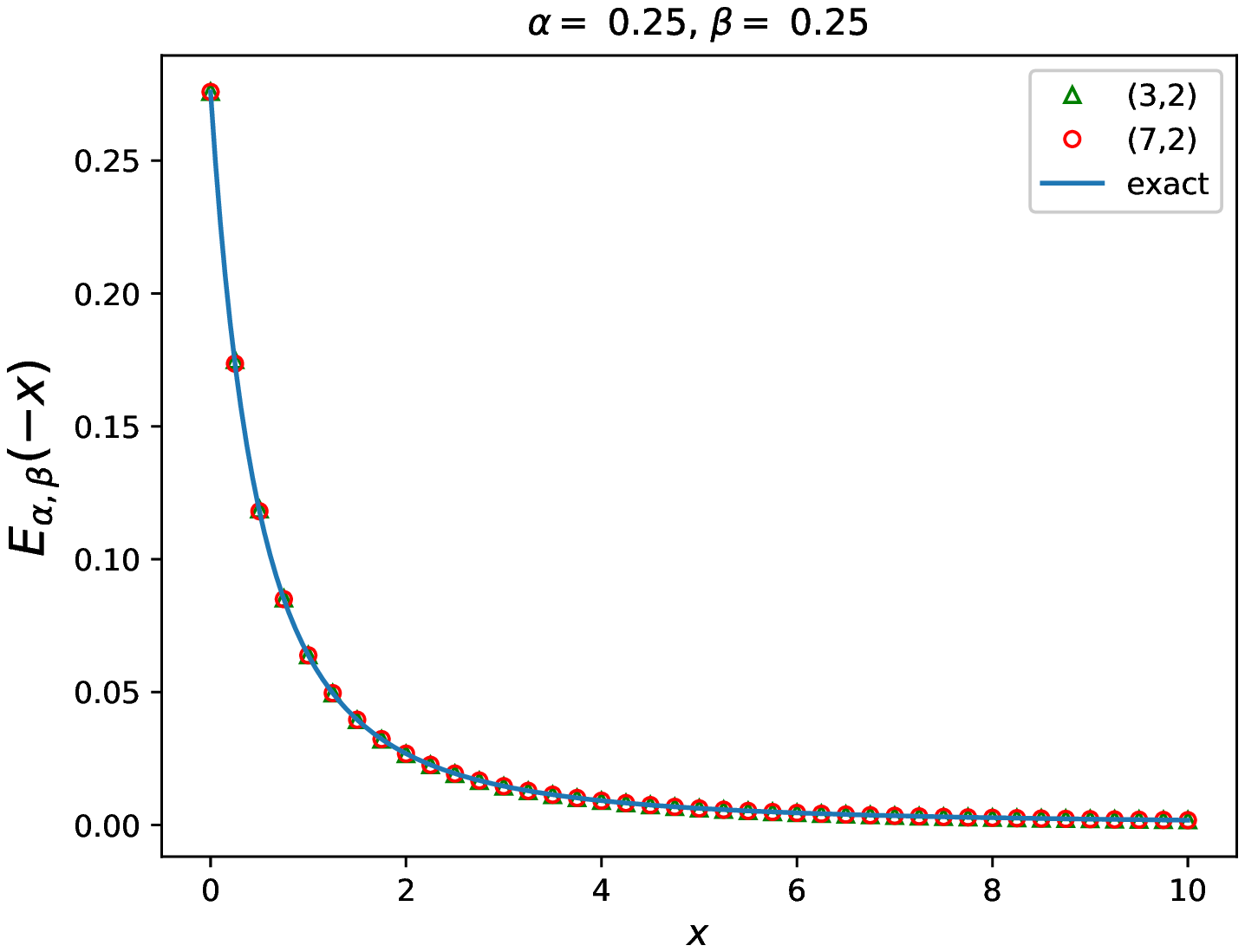}
\includegraphics[width=0.49\textwidth]{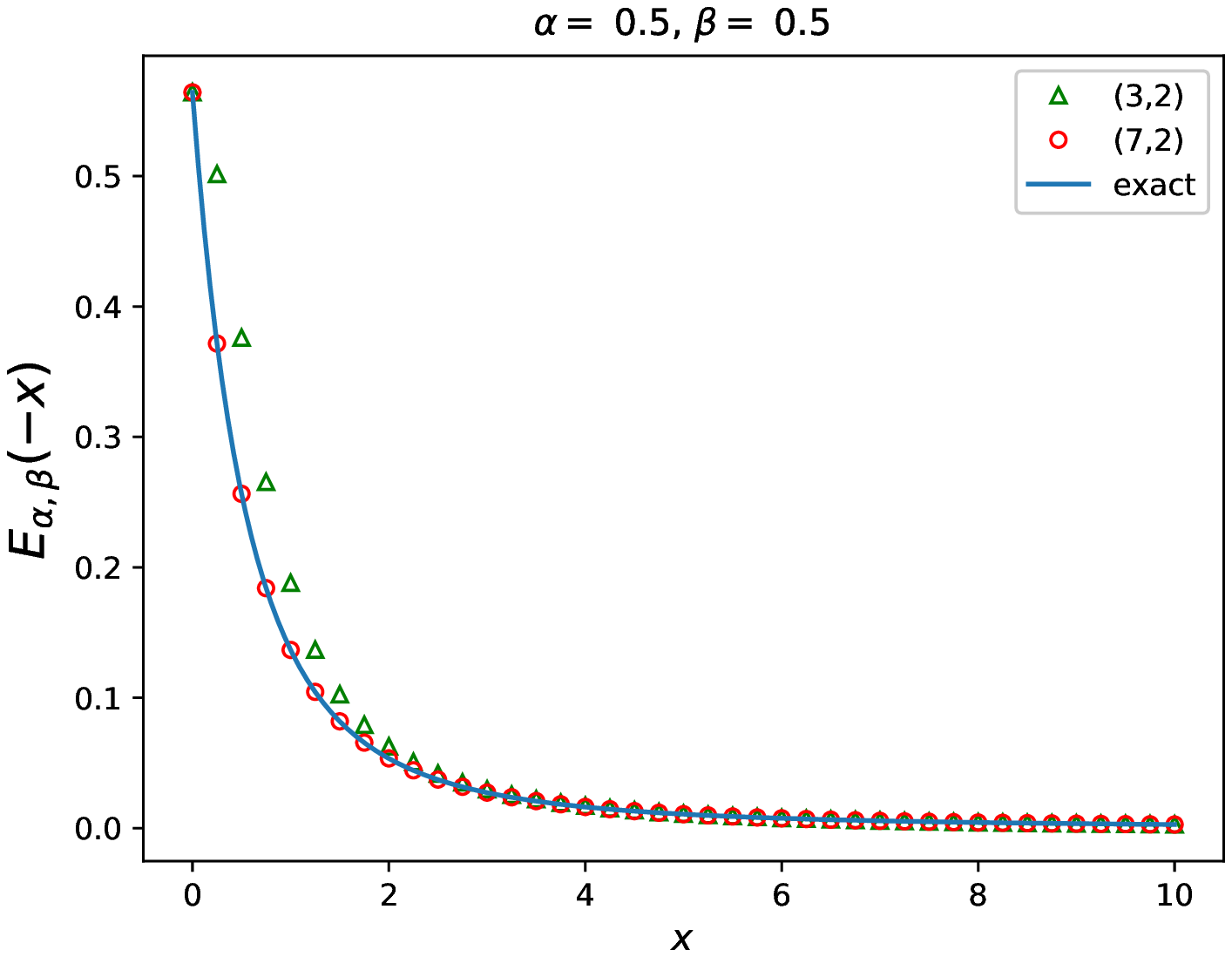} 
\\
\includegraphics[width=0.49\textwidth]{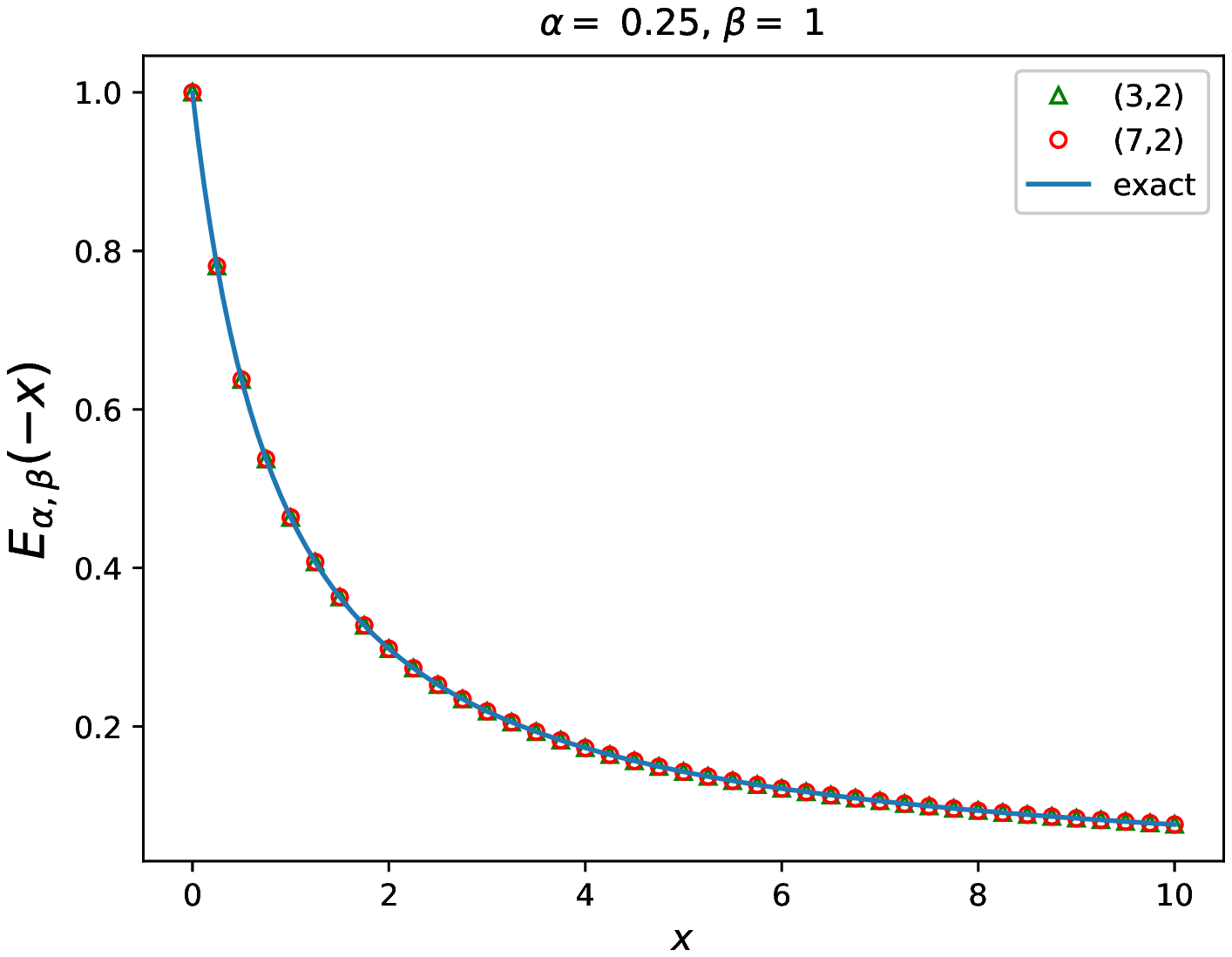}
\includegraphics[width=0.49\textwidth]{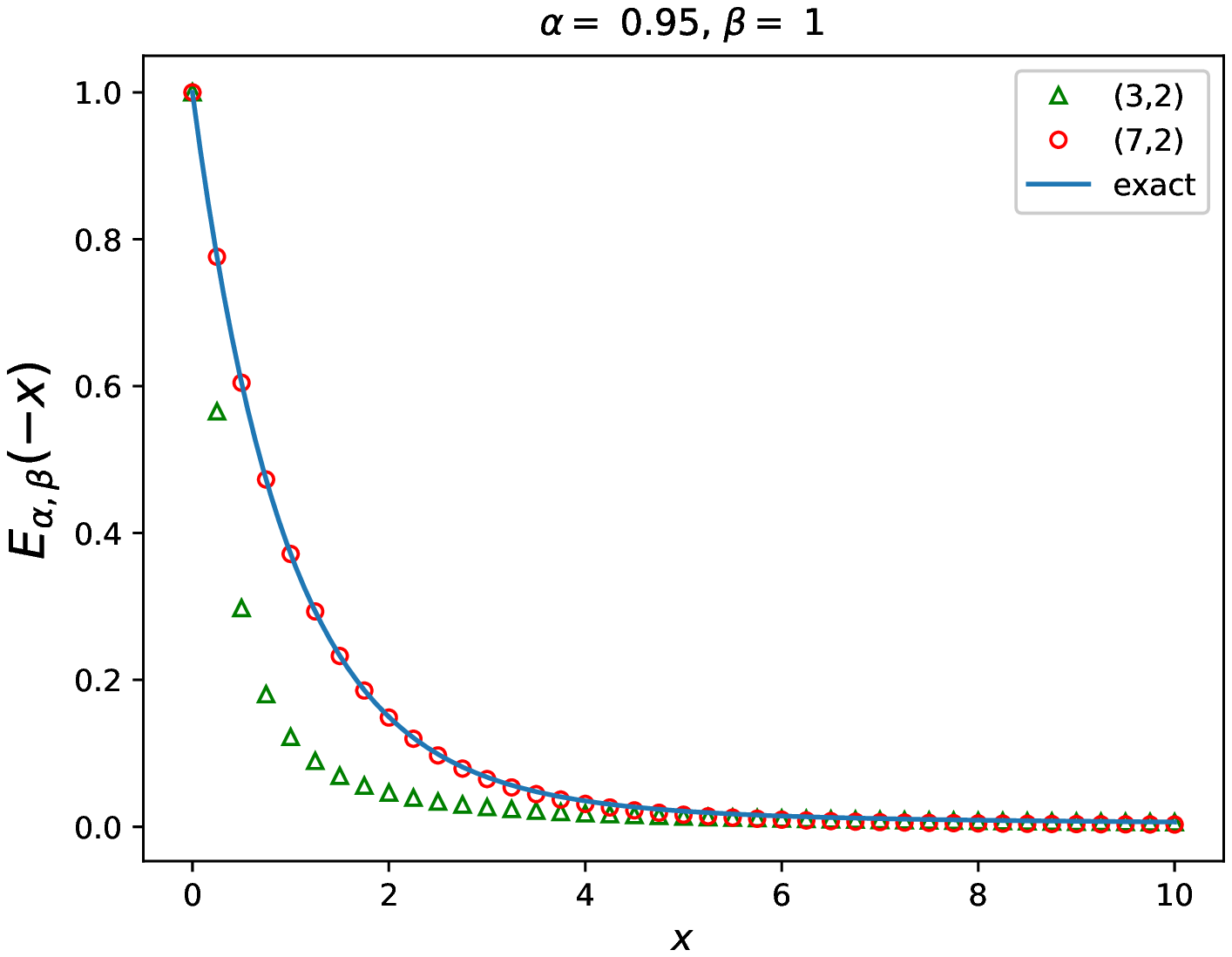} 
\\
\includegraphics[width=0.49\textwidth]{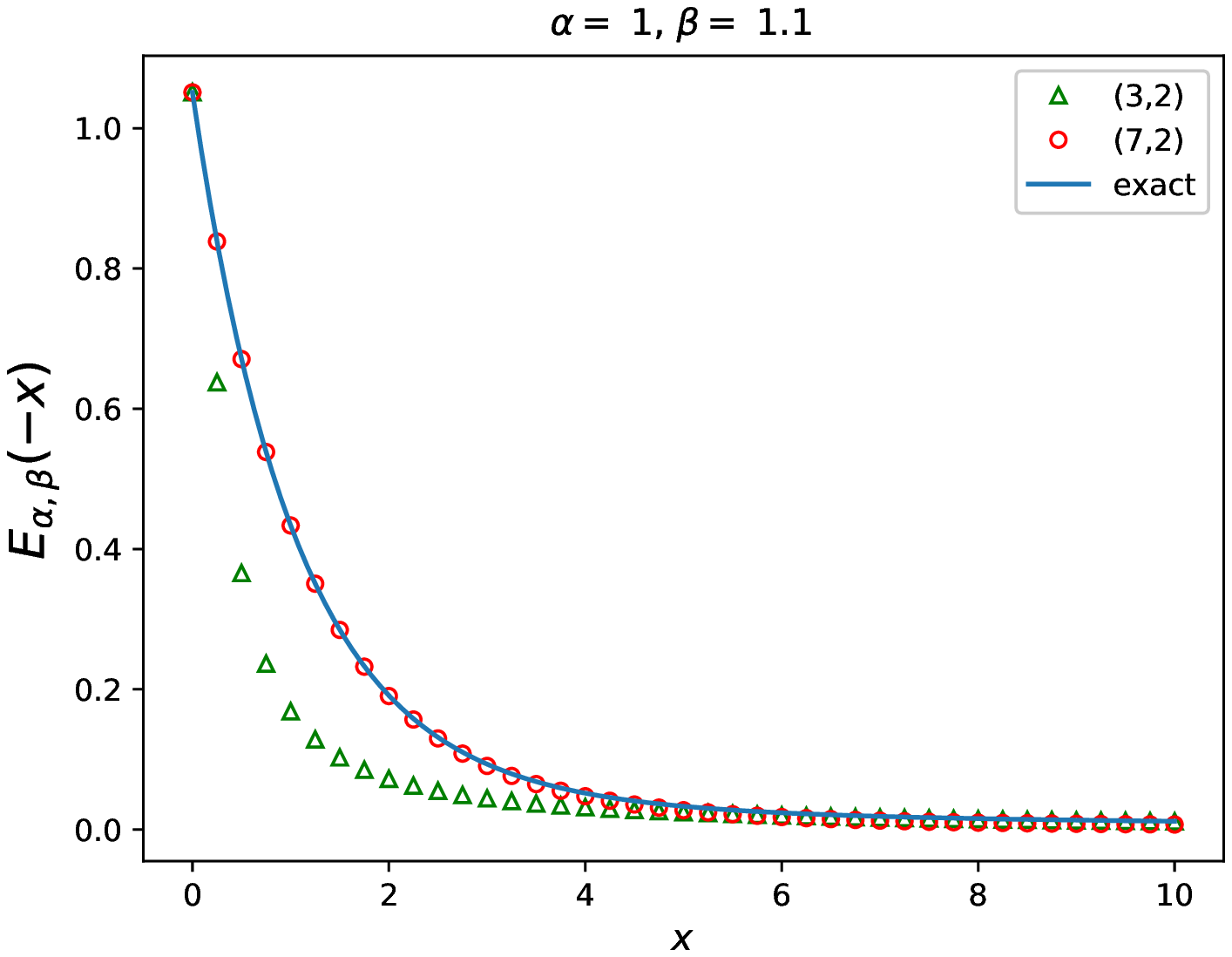}
\includegraphics[width=0.49\textwidth]{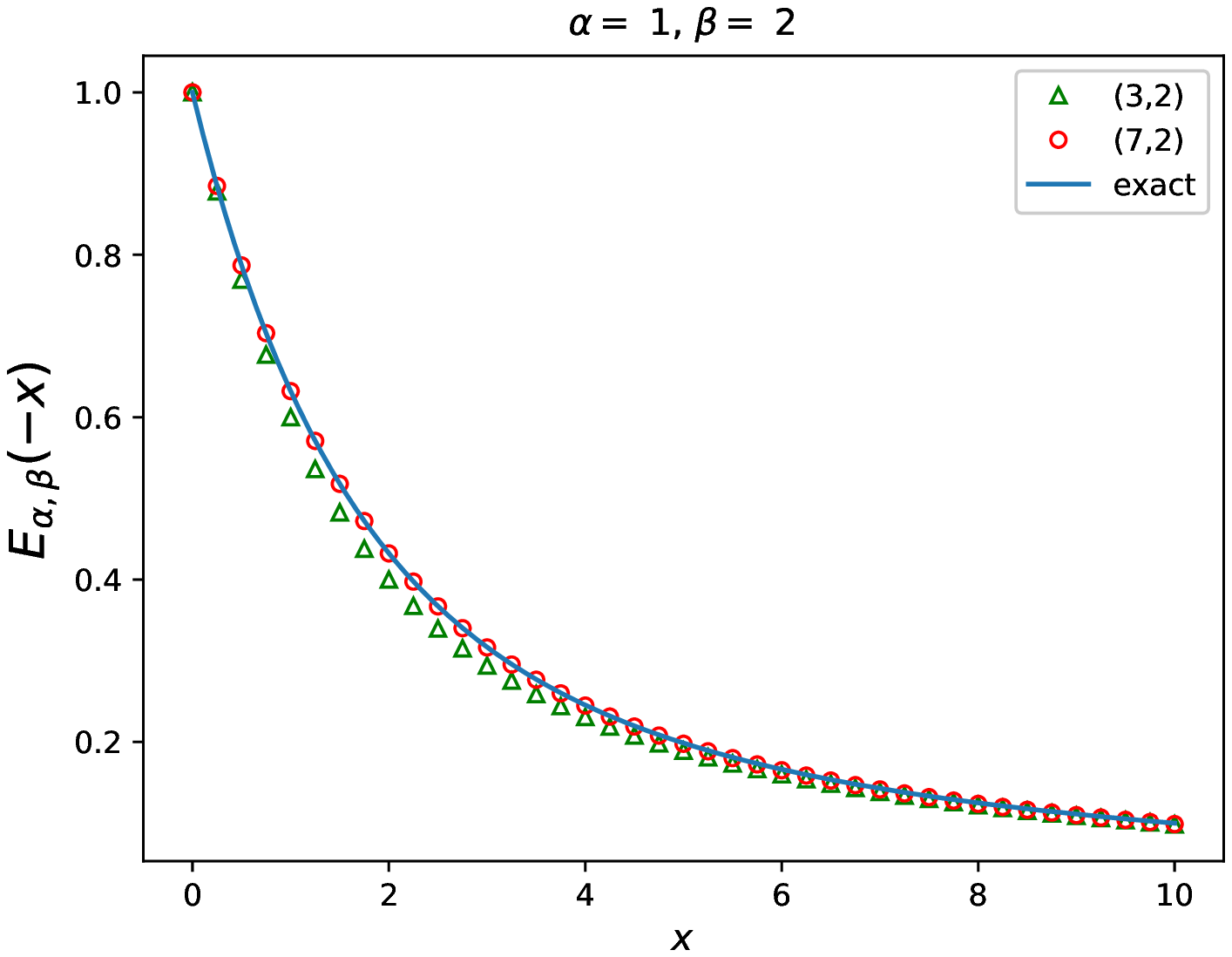}
\caption{Comparisons of $R_{\alpha,\beta}^{7,2}$ and $R_{\alpha,\beta}^{3,2}$}
\label{fig:fourth order vs second order}
\end{figure}

In this section, we demonstrate graphically and computationally the performance of the fourth-order global Pad\'e approximants constructed above.
We start by comparing the approximants $R_{\alpha, \beta}^{5,4}$,  $R_{\alpha, \beta}^{6,3}$ and $R_{\alpha, \beta}^{7,2}$. 
Figures \ref{fig:fourth order beta 1}, \ref{fig:fourth order alpha 1} and 
\ref{fig:fourth order alpha is beta} contain the profiles for different combinations of $(\alpha,\beta)$.  
These profiles reveal that both $R_{\alpha, \beta}^{7,2}(x)$ and $R_{\alpha, \beta}^{6,3}(x)$ 
provide extremely well approximants and compare favorably with
$R_{\alpha, \beta}^{5,4}(x)$.
This observation supports our earlier argument in Remark \ref{remark:error} that approximations of the same order can be improved by increasing the 
number of local terms $m$ and decreasing the number of asymptotic terms $n$ in the matching requirements \eqref{eq:gpa matching} and 
\eqref{eq:gpa matching2}.
The profile of the absolute errors are shown in Figure \ref{fig:fourth order approximation errors} and for completeness, we include in 
Figure~\ref{fig:fourth order vs second order} different profiles of $R_{\alpha,\beta}^{7,2}$ vs $R_{\alpha,\beta}^{3,2}$.

In table \ref{tab:error table}, we provide the maximum absolute error 
$$
\max\limits_{x \in I}\{\abs{E_{\alpha, \beta}(-x) - R_{\alpha, \beta}(x)}\},
$$
and the maximum relative error
$$
\max\limits_{x \in I}\bigg\{\abs{\dfrac{E_{\alpha, \beta}(-x) - R_{\alpha, \beta}(x)}{E_{\alpha, \beta}(-x)}}\bigg\}, 
$$
where $I$ is the interval $[0,10]$ and $R_{\alpha, \beta}$ is any of the approximants for $E_{\alpha,\beta}$. 
In the case $\beta = 1$, the errors resulting from $R_{\alpha, 1}^{6,3}$ and 
$R_{\alpha, 1}^{7,2}$ are smaller than those from the approximation by Ingo \etal in \cite{Ingo:2017}. 
It is worth mentioning that they were comparing with "mlf" Matlab function by Podlubny which is based on the algorithm by Gorenflo \etal in \cite{Gorenflo:2002} while we are comparing with the "ml" Matlab function by Garrappa \cite{Garrappa:2015}.

\begin{figure}
\centering
\includegraphics[width=0.49\textwidth]{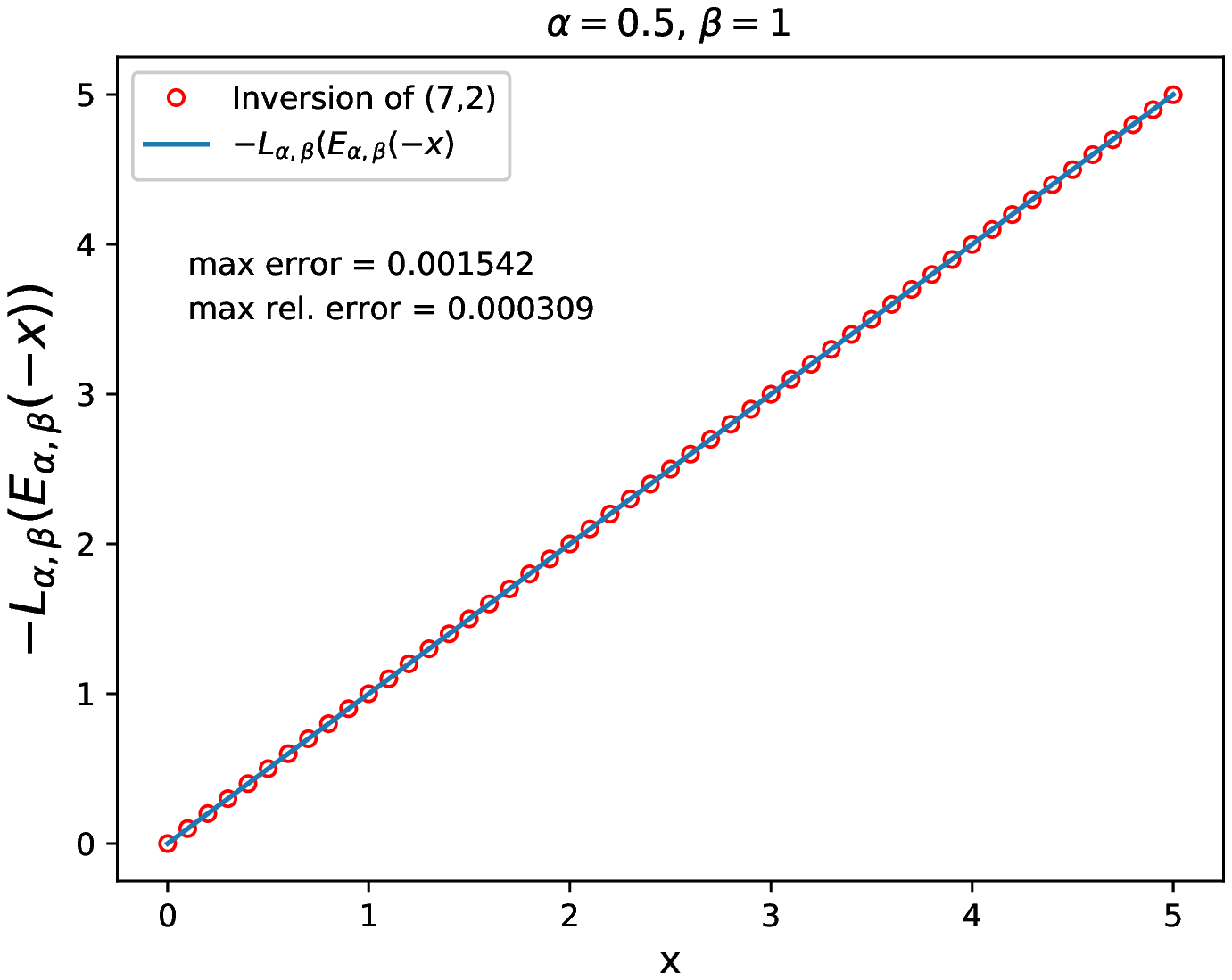}
\includegraphics[width=0.49\textwidth]{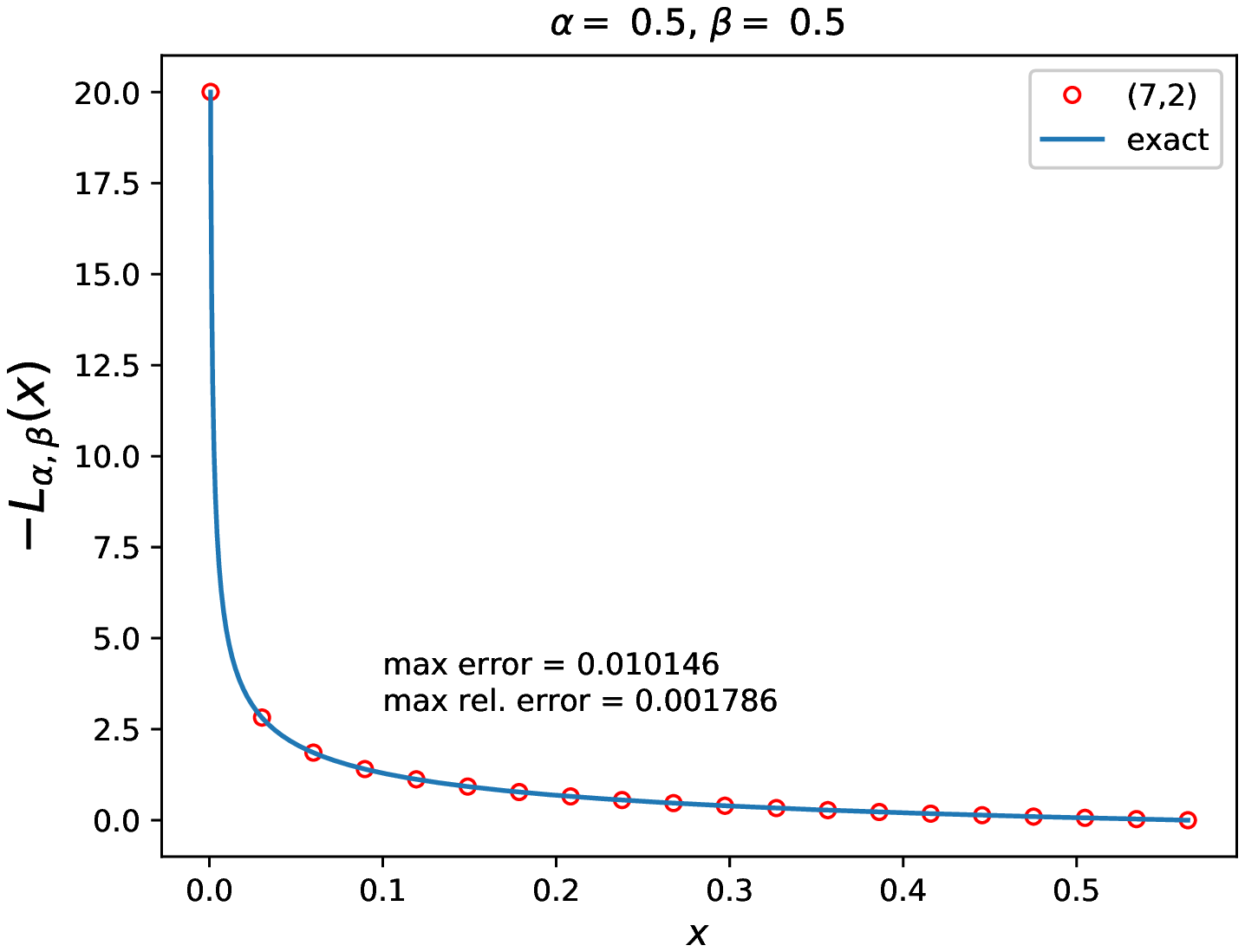}
\\
\includegraphics[width=0.49\textwidth]{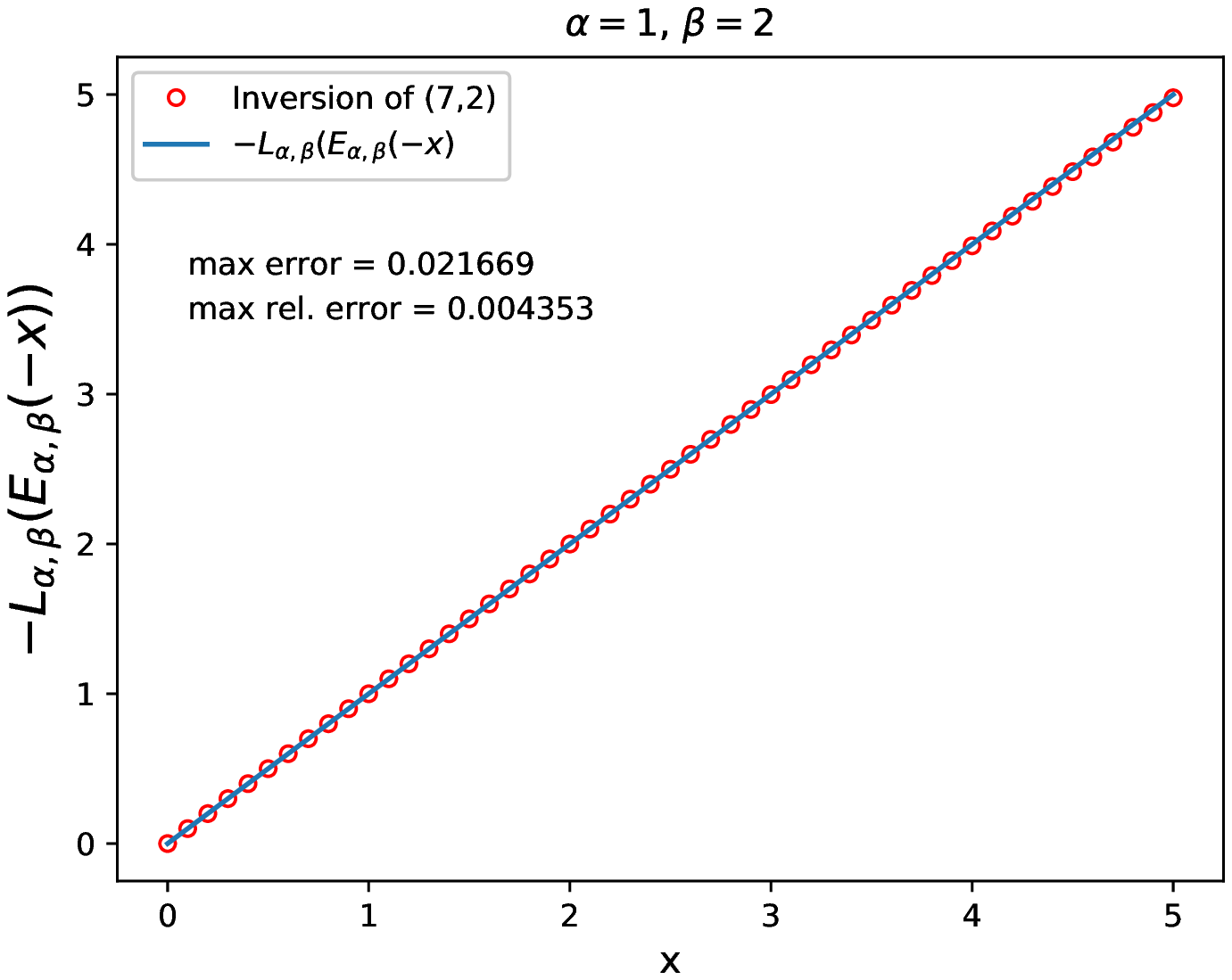} 
\includegraphics[width=0.49\textwidth]{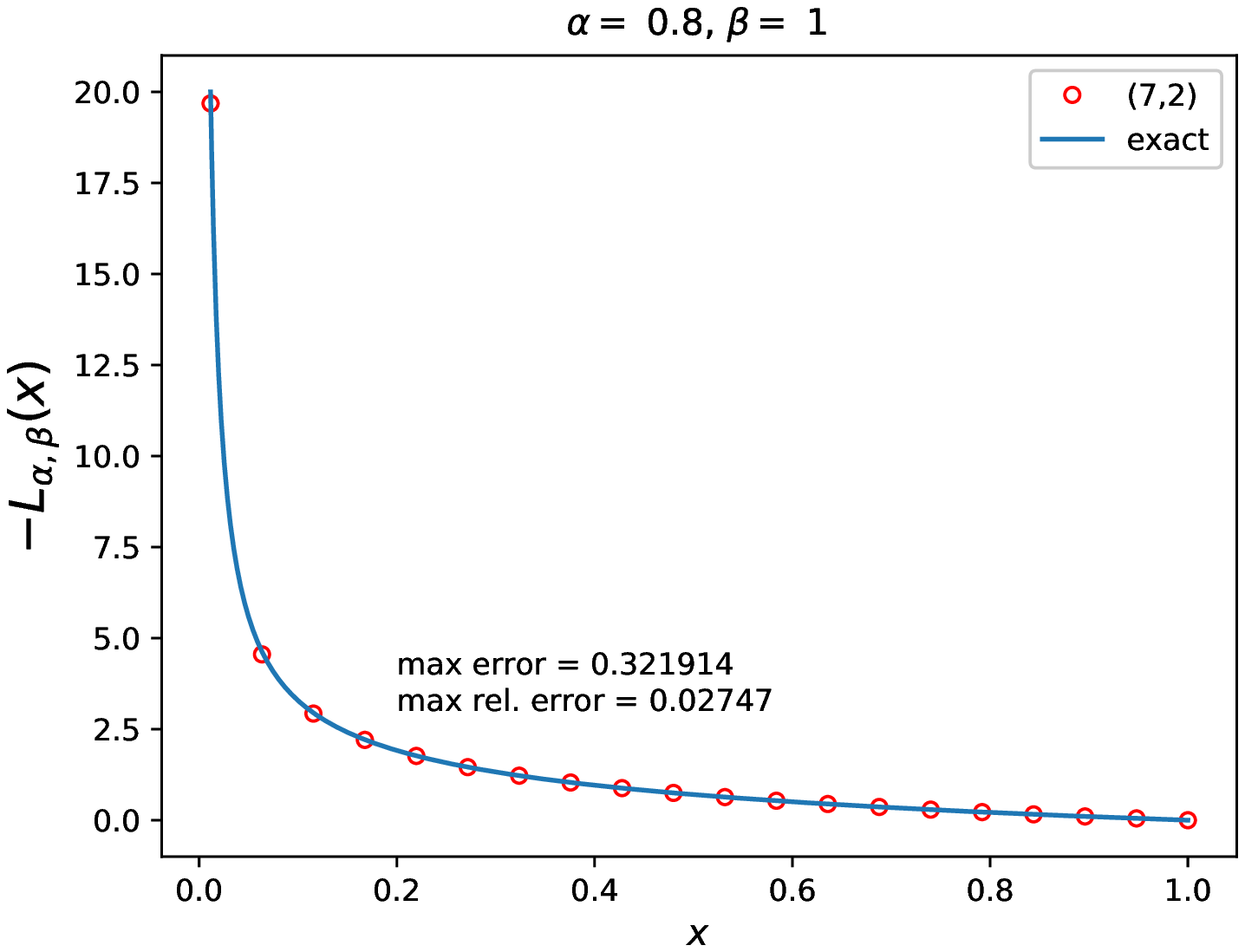} 
\caption{Plots of inverse MLF approximated by solving \eqref{eq:inverse gpa} 
for $R^{7,2}_{\alpha, \beta}$ vs the exact value obtained by applying the definition \eqref{eq:inverse mlf}.}
\label{fig:inverse mlf}
\end{figure}

Finally, in Figure \ref{fig:inverse mlf} we compare the approximations of the inverse function obtained by solving \eqref{eq:inverse gpa} 
for $R^{7,2}_{\alpha, \beta}$ vs the exact values by definition \eqref{eq:inverse mlf}.

\begin{table}[tb]
\begin{center}
\small{
\begin{tabular}{c| cc|cc|cc|cc|cc}
&\multicolumn{2}{c|}{$\alpha = 0.9, \beta = 1.9$} & \multicolumn{2}{c|}{$\alpha = 0.9, \beta = 1$} & 
\multicolumn{2}{c|}{$\alpha = 0.5, \beta = 0.5$} &
\multicolumn{2}{c|}{$\alpha = 1.0, \beta = 1.1$} &
\multicolumn{2}{c}{$\alpha = 1.0, \beta = 2.0$}
\\
& AE  & RE &  AE  & RE   & AE  &   RE    & AE &   RE  & AE & RE  
\\ \hline 
& & & & & & & & & \\
RDP	& 0.0813  & 0.8234   & 0.1153  & 3.2415   & 0.0800  &  10.675   
    & 0.1225  &  3.3583  & 0.1165  & 1.165  
\\ & & & & & & & & & \\
$R_{\alpha,\beta}^{3,2}$ & 0.0254  & 0.0526   & 0.1926  & 0.4884   
& 0.1349  &  0.4697      & 0.3112  &  0.6394  & 0.0352  & 0.0755
\\ & & & & & & & & & \\
$R_{\alpha,\beta}^{5,4}$ & 0.0018   & 0.0052   & 0.0264  & 0.1552  
&  0.0040  &  0.0119     & 0.1020   &  0.4274  & 0.0040  & 0.0119
\\ & & & & & & & & & \\
$R_{\alpha,\beta}^{6,3}$ & 0.0007 & 0.0034  & 0.0047  & 0.0920   
& 0.0002  &  0.0028      & 0.0092 & 0.2278  & 0.0014  & 0.0073 
\\ & & & & & & & & & \\
$R_{\alpha,\beta}^{7,2}$ & 0.0006 & 0.0055   & 0.0033  & 0.1648   
& 0.0001  &  0.0033      & 0.0061 &  0.3949  & 0.0012  & 0.0112   
\\ 
\end{tabular} 
}
\caption {The maximum absolute error (AE) and maximum relative error (RE) for the different approximants}
\label{tab:error table}
\end{center}
\end{table} 

\newpage
\section{Applications}
\label{sec:applcications}
The fundamental importance of the two-parametric MLF and its inverse is the 
main motivation behind the construction of the approximants in this paper. 
Our objective here is to show the accuracy and efficiency of these approximants when applied to solutions of fractional differential 
and integral equations. 
The fourth order approximant $R^{7,2}_{\alpha, \beta}$ 
is used all through this section, while the second order approximant $R^{3,2}_{\alpha, \beta}$ is also applied for comparison purposes.
Below, we consider some applications with solutions that involve MLF with scalar and matrix arguments.

\subsection{Applications with scalar arguments}

We start by considering the following applications that involve MLF with scalar arguments.

\subsubsection{Fractional reaction-diffusion equation}

\begin{figure}
\centering
\includegraphics[width=0.49\textwidth]{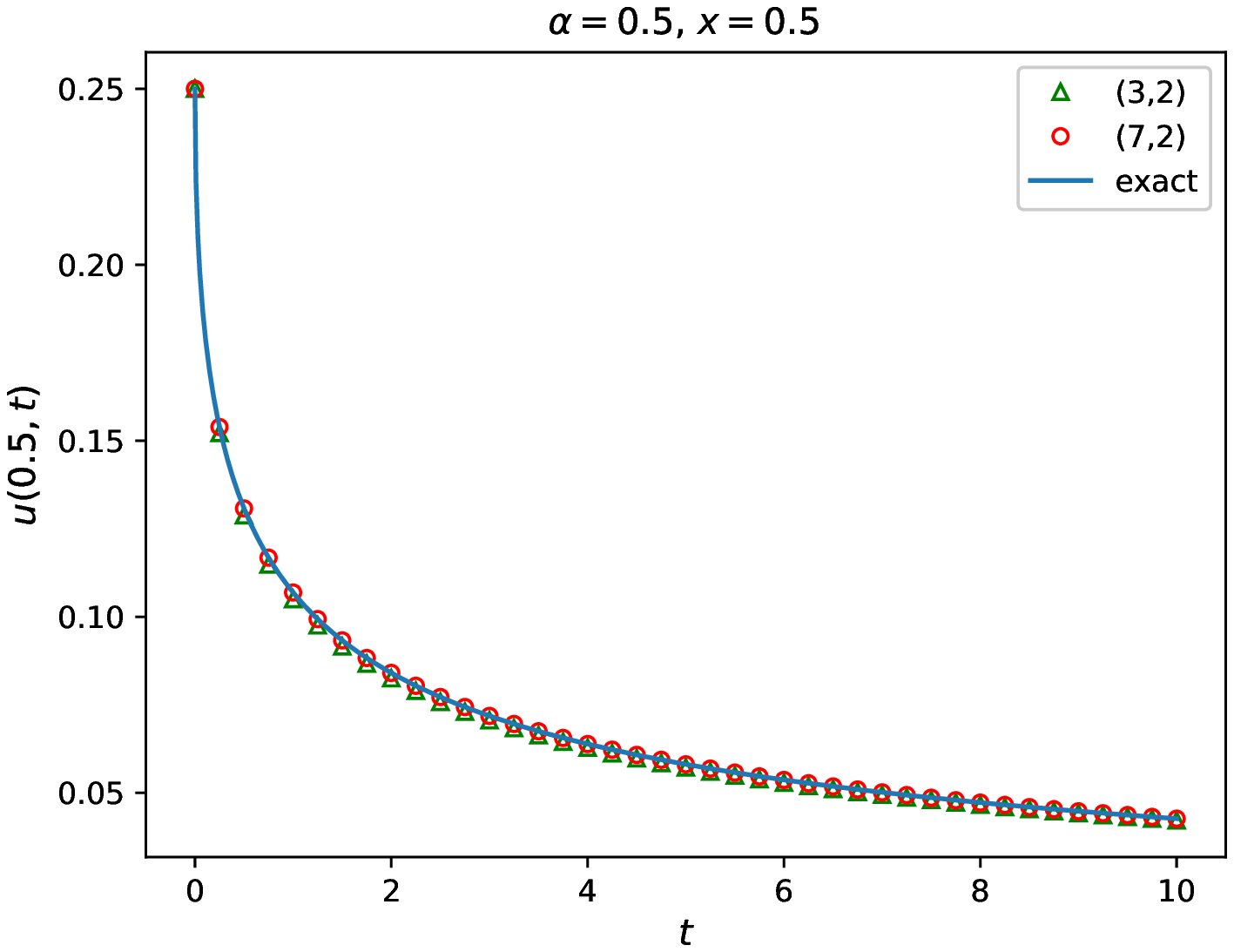}
\includegraphics[width=0.49\textwidth]{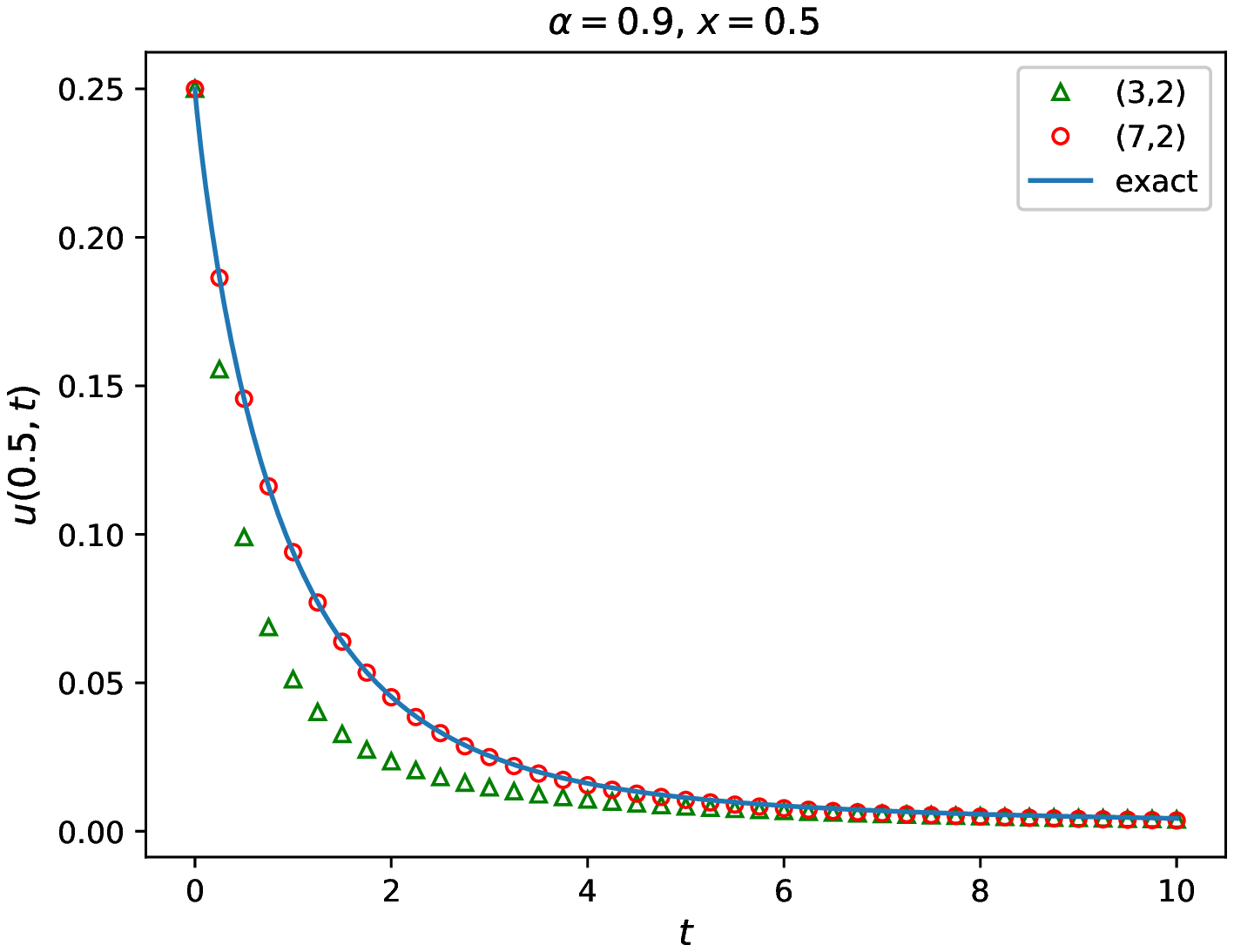}
\caption{
Plots of $R_{\alpha, \beta}^{7,2}$ and $R_{\alpha, \beta}^{3,2}$ approximations to the solution \eqref{eq:reaction diffusion solution} 
of the fractional reaction-diffusion problem~\eqref{eq:reaction diffusion} at $x=0.5$
}
\label{fig:reaction-diffusion}
\end{figure}

Consider the following sub-diffusion initial-boundary value problem:
\begin{align}
\label{eq:reaction diffusion}
& \cD^\alpha_t u(x,t) = u_{xx} (x,t) + u(x,t) \, u_x (x,t) + f(x,t), 
\quad x\in (0,1), t > 0,
\\ & \nonumber
u(0,t) = u(1,t) = 0, \quad t \geq 0,
\\ & \nonumber
u(x,0) = x(1-x), \qquad x \in [0,1],
\end{align}
where $\cD^\alpha_t$, $0<\alpha<1$, is the Caputo fractional derivative.
When
$$
f(x,t) = - \Big[ 2 + E_{\alpha}(- t^\alpha)(2x^3 - 
3x^2 + x) + x(1-x) \Big] \, E_{\alpha}(- t^\alpha),
$$ 
then the exact solution is 
\begin{equation}
\label{eq:reaction diffusion solution}
u(x,t) = x(1-x) E_{\alpha}(- t^\alpha).
\end{equation}
The solution profile and its approximations at $x=0.5$ for $t\in [0,10]$ and $\alpha = 0.5, 0.9$ are included in Figure \ref{fig:reaction-diffusion}.
The corresponding errors and runtime (in seconds) for time increment of $0.01$
are listed in Table~\ref{tab:comparison RDE} together with the runtime for ml subroutine in \cite{Garrappa:2018}.
As can be observed, $R_{\alpha, \beta}^{7,2}$ provides an excellent approximation at a significantly reduced runtime as compared to the ml subroutine runtime.
Furthermore, it is clear from the figures that, in general, $R_{\alpha,\beta}^{3,2}$ may not be a good option.

\begin{table}[h]
\centering
\begin{tabular}{l| ccc|c} 
	$\alpha$  & Max Abs. Error  & Max Rel. Error & Runtime & ml Runtime
	\\ \hline & & & & \\
	$0.5$ & 9.56e-06 & 6.77e-05 & 6.22e-04  & 7.25e-02 
	\\
	$0.9$ & 8.25e-04 & 1.80e-03 & 7.06e-04   & 9.66e-02  
	\\  & & & & \\ \hline 
\end{tabular} 
\caption{
Errors and runtime for $R_{\alpha, \beta}^{7,2}$ approximation to the solution 
\eqref{eq:reaction diffusion solution} of the fractional reaction-diffusion problem \eqref{eq:reaction diffusion} at $x=0.5$
}
\label{tab:comparison RDE}
\end{table} 

\subsubsection{Fractional integral equation}
\label{VIE}

\begin{figure}
\includegraphics[width=0.49\textwidth]{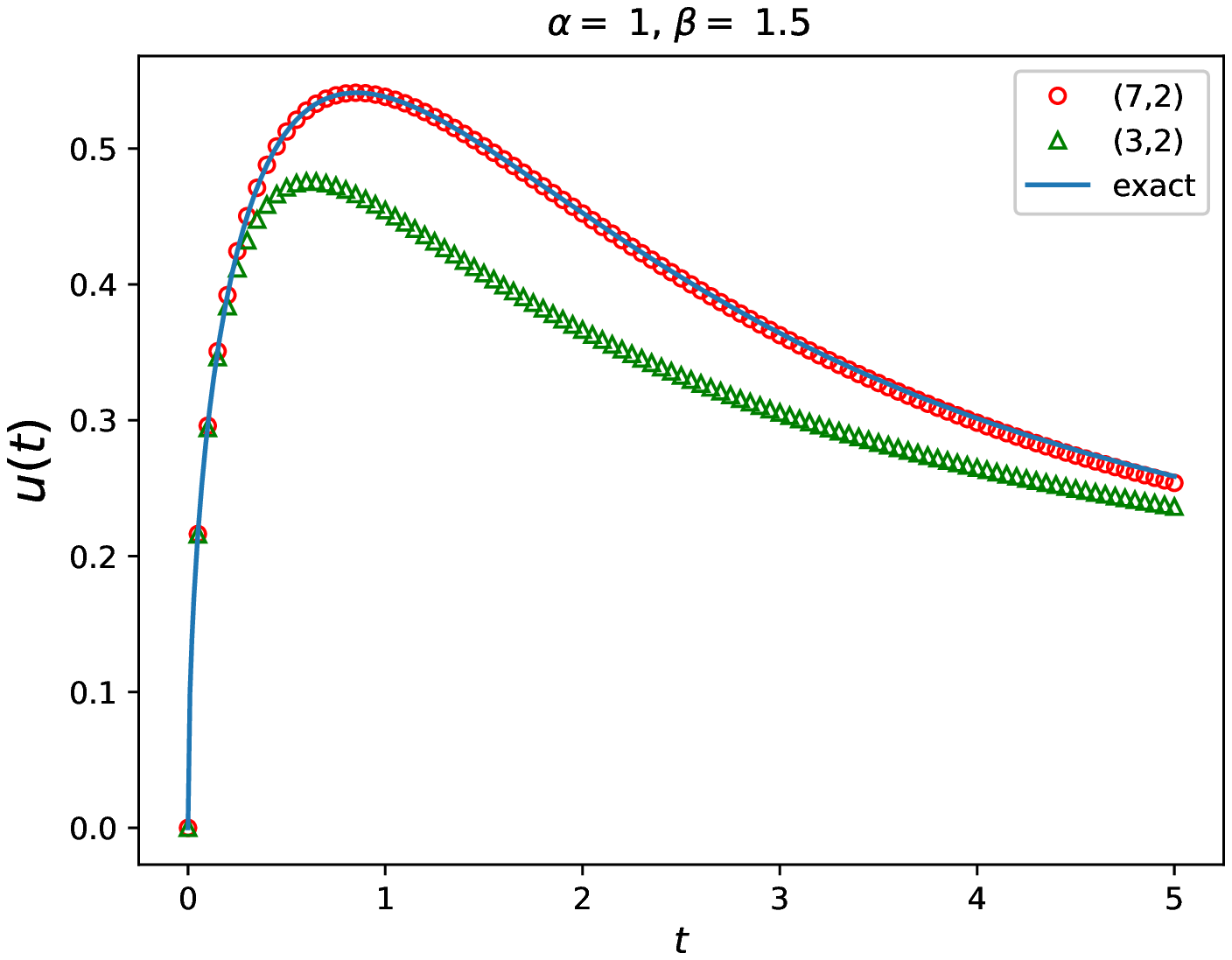}
\includegraphics[width=0.49\textwidth]{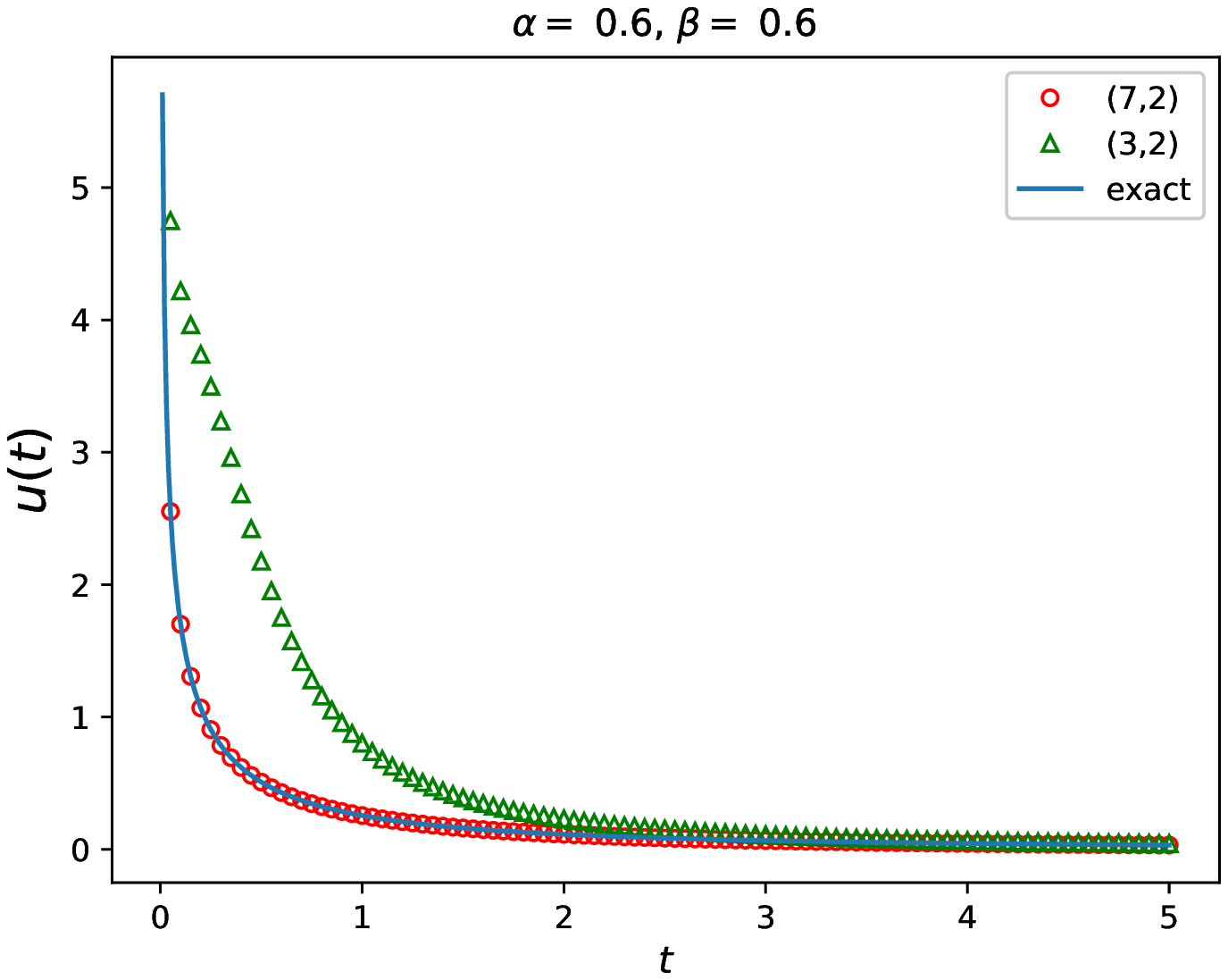}\\
\caption{
Plots of $R_{\alpha, \beta}^{7,2}$ and $R_{\alpha, \beta}^{3,2}$ approximations to the solution 
\eqref{eq:VIE solution} of the integral equation~\eqref{eq:VIE}
}
\label{fig:VIE}
\end{figure}

Consider the following integral equation:
\begin{equation}\label{eq:VIE}
\mathcal{I}^\alpha_t u(t) + u(t) =  t^{\beta - 1}, \quad t > 0,
\end{equation}
where $\alpha \text{, } \beta > 0$ and $\mathcal{I}^\alpha$ is the Riemann-Liouville integral of order $\alpha$. 
The exact solution of \eqref{eq:VIE} is given by:
\begin{equation}\label{eq:VIE solution}
u(t) = \Gamma(\beta) t^{\beta - 1} E_{\alpha, \beta} (-t^\alpha).
\end{equation} 
Plots of the solution $u(t)$ and its approximations for 
$(\alpha,\beta) = (0.6, 0.6)$ and $(1.0, 1.5)$ 
are provided in Figure \ref{fig:VIE}.
The corresponding errors and runtime (in seconds) for time increment of $0.01$ are listed in Table~\ref{tab:comparison VIE}.
Again, the results assert that superiority of $R_{\alpha, \beta}^{7,2}$.

\begin{table}[h]
\centering
\begin{tabular}{l| ccc|c} 
$(\alpha,\beta)$  & Max Abs. Error  & Max Rel. Error & Runtime & ml Runtime
\\ \hline & & & & \\
$(0.6, 0.6)$ & 3.77e-04 & 1.39e-04 & 3.00e-03  & 8.09e-02 
\\
$(1.0,1.5)$ & 7.40e-03 & 8.70e-03 & 1.20e-03   & 1.09e-01  
\\  & & & & \\ \hline 
\end{tabular} 
\caption{
Errors and runtime for $R_{\alpha, \beta}^{7,2}$ approximation to the solution 
\eqref{eq:VIE solution} of the integral equation \eqref{eq:VIE}
}
\label{tab:comparison VIE}
\end{table} 

\subsubsection{Ultraslow diffusion}
 
The propagator $p(x,t)$ of an ultraslow diffusive process satisfies the structural diffusion equation:
\begin{equation}
\label{eq:ultraslow diffusion}
\frac{dp(x,t)}{d_mt} = k_\alpha \partial^2_xp(x,t), \qquad t > 0, -\infty < x < \infty,
\end{equation}
where the local structural derivative in time $\frac{dp(x,t)}{d_mt}$ 
with respect to the structural function $E_{\alpha}^{-1}$ is given by \cite{Liang:2018}
$$
\frac{dp(x,t)}{d_mt} = 
\lim_{s \to t} \frac{p(x,s) - p(x,t)}{E_\alpha^{-1}(s) - E_\alpha^{-1}(t)},
\qquad 0 < \alpha < 1.
$$
The solution of \eqref{eq:ultraslow diffusion} (see \cite{Liang:2018}) is given by the scaled Gaussian function:
\begin{equation}
\label{eq:ultraslow solution}
p(x,t) = \frac{1}{\sqrt{4\pi k_\alpha E^{-1}_\alpha(t)}}\exp(-\frac{x^2}{4 k_\alpha E^{-1}_\alpha(t)}).
\end{equation}

\begin{figure}
\centering
\includegraphics[width=0.49\textwidth]{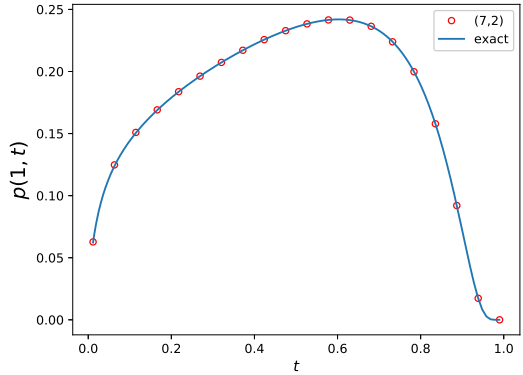}
\caption{
Plots of the propagator $p(x,t)$ \eqref{eq:ultraslow solution}, of the ultraslow diffusion equation \eqref{eq:ultraslow diffusion}
together with the corresponding approximation using the inversion of $R^{7,2}_{\alpha, \beta}$ at $x = 1$ with $\alpha = 0.6$.
}
\label{fig:ultraslow}
\end{figure}

As shown in Figure \ref{fig:ultraslow}, we have a good agreement between the propagator 
$p(x,t)$ for $t\in(0,1)$, $x = 1$, and its approximation when $E^{-1}_\alpha$
is approximated by the inverse of $R_{\alpha, 1}^{7,2}$. 

\subsection{Applications with matrix arguments}
MLFs of matrix argument arise naturally in solutions of systems of fractional 
differential equations. 
The two-parametric MLF of a matrix 
$A \in \mathbb{C}^{n \times n}$ (see \cite{Sadeghi:2018}) is given by:
\begin{equation}
\label{eq:mlf matrix}
E_{\alpha, \beta}(A) = \sum_{k=0}^{\infty} \frac{A^k}{\Gamma(\alpha k + \beta)}, 
\qquad \Re \alpha > 0, \quad \beta \in \C.
\end{equation}

Next, we demonstrate the effectiveness of using 
$R^{7,2}_{\alpha, \beta}$ to approximate $E_{\alpha,\beta}$ of a matrix argument. 
The numerical experiments below show the runtime saving and accuracy 
in comparison with the ml\_matrix function described in \cite{Garrappa:2018}.

\subsubsection{The Bagley-Torvik problem}

The initial-value problem of Bagley-Torvik equation is given by
\begin{equation}\label{eq:Bagley-Torvik}
D^2 u(t) + a_1 \OcD^{3/2} u(t) + a_2 u(t) = f(t), 
\quad u(0) = 0, \, u'(0) = 0,
\end{equation}
where $a_1$ and $a_2$ are positive constants.
This problem models the motion of a thin, rigid plate 
immersed in a Newtonian fluid of infinite extension connected 
to a fixed point via a spring \cite{Kilbas:2006b}.
Using Laplace transform, the exact solution of \eqref{eq:Bagley-Torvik} is
\begin{equation}
u(t) = \int_{0}^{t}G_{2,\frac{3}{2}}(t-\tau) f(\tau) d\tau,
\end{equation}
where 
\begin{equation}\label{eq:hypergeom}
G_{2,\frac{3}{2}}(t) = \sum_{n = 0}^{\infty} \frac{(-a_2)^n}{n!} t^{2n+1} 
\wfun{1}{1} \left[ \left. \begin{array}{c} (n+1, 1) \\ \\(2n+2, \frac{1}{2})
\end{array} \right| -a_1 t^\frac{1}{2} \right].
\end{equation}
and $\wfun{1}{1}$ is the Wright function.

To avoid the complexity of calculating \eqref{eq:hypergeom}, 
using Theorem 8.1 in \cite{Diethelm:2010b}, the initial-value problem \eqref{eq:Bagley-Torvik} 
can be converted into the system
\begin{equation}
\label{eq:Bagley-Torvik system}
D^{\frac{1}{2}}U(t) = AU(t) + f(t) \, e_4, \qquad U(0) = U^0,
\end{equation}
where
$$
U = \left(u, \OcD^\frac{1}{2} u, \OcD^1 u, \OcD^\frac{3}{2} u \right)^T, \quad
e_4 = (0, 0, 0, 1)^T, \qquad U^0 = (0, 0, 0, 0)^T,
$$ 
and
$$ 
A = \begin{bmatrix}
0  & 1  & 0 &  0     \\
0  & 0  & 1 &  0     \\
0  & 0  & 0 &  1      \\
-a_2 & 0  & 0 & -a_1
\end{bmatrix}.
$$
Then, the exact solution of \eqref{eq:Bagley-Torvik system} is \cite{Popolizio:2018}
\begin{equation}
\label{eq: system solution}
U(t) = \int_0^t (t-\tau)^{-\frac{1}{2}} \, E_{\frac{1}{2},\frac{1}{2}} ((t-\tau)^\frac{1}{2} A) \, e_4 \, f(\tau) \, d\tau.
\end{equation}

For testing purposes, let
\begin{equation}\label{eq:Bagley-Torvik source}
f(t) = a_2t^2 + \frac{4 a_1}{\sqrt{\pi}} t^\frac{1}{2} + 2.
\end{equation}
Then the exact solution of \eqref{eq:Bagley-Torvik} is 
\begin{equation}
\label{eq:Bagley-Torvik solution}
u(t) = t^2,
\end{equation}
while the solution \eqref{eq: system solution} takes the form
\begin{equation}
\label{eq:Bagley-Torvik system solution}
U(t) = 2 \left[  t^\frac{1}{2} E_{\frac{1}{2},\frac{3}{2}}(At^\frac{1}{2}) + 
a_1 t E_{\frac{1}{2},2}(At^\frac{1}{2}) +  
a_2 t^\frac{5}{2} E_{\frac{1}{2},\frac{7}{2}}(At^\frac{1}{2}) \right]  e_4.
\end{equation}
Then the MLFs of matrix argument in \eqref{eq:Bagley-Torvik system solution} can be approximated using the 
partial fraction decomposition \eqref{eq:partial-frac conjugate fourth order} of $R_{\alpha,\beta}^{7,2}$:
$$
E_{\alpha,\beta}(B) \approx 2\Re\left[c_1(-B-Ir_1)^{-1} + c_2(-B-Ir_2)^{-1} \right].
$$
Explicitly, we have
\begin{align*}
2 \, t^\frac{1}{2} E_{\frac{1}{2},\frac{3}{2}}(At^\frac{1}{2}) \, e_4 
& \approx 
4 \Re \bigg( t^\frac{1}{2}\left[c_1(-At^\frac{1}{2}-Ir_1)^{-1} \right] e_4 + t^\frac{1}{2}\left[c_2(-At^\frac{1}{2}-Ir_2)^{-1}\right] e_4\bigg)
\\ & =
4\Re [v_1 + v_2 ],
\end{align*}
where the vectors $v_1$ and $v_2$ are obtained by solving the systems
$$
(-At^\frac{1}{2} - Ir_1) \, v_1 = c_1 t^\frac{1}{2} e_4, \qquad
(-At^\frac{1}{2} - Ir_2) \, v_2 = c_2 t^\frac{1}{2} e_4.
$$
The remaining terms in \eqref{eq:Bagley-Torvik system solution} are approximated in a similar manner.

\begin{figure}
\includegraphics[width=0.49\textwidth]{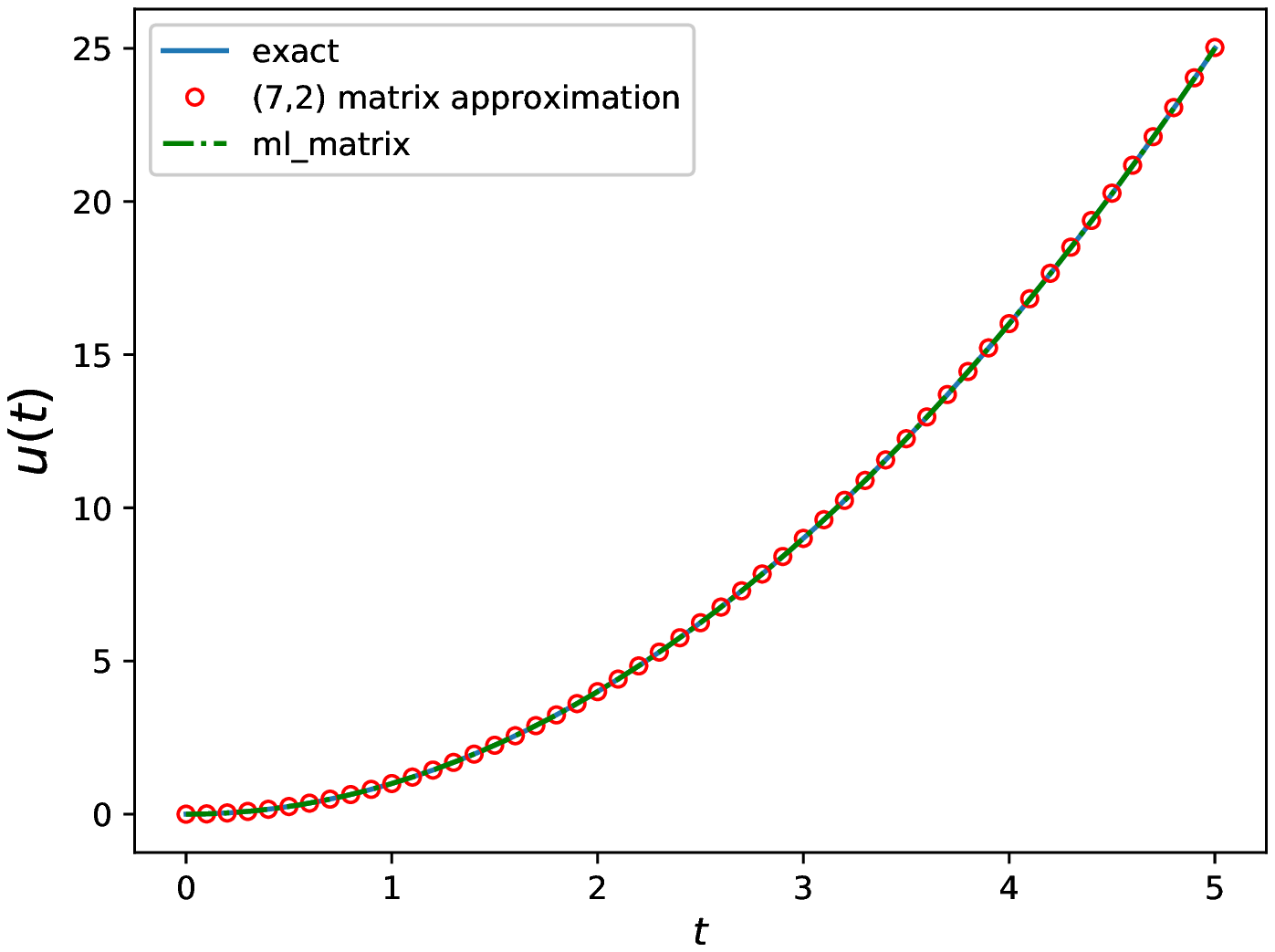}
\includegraphics[width=0.49\textwidth]{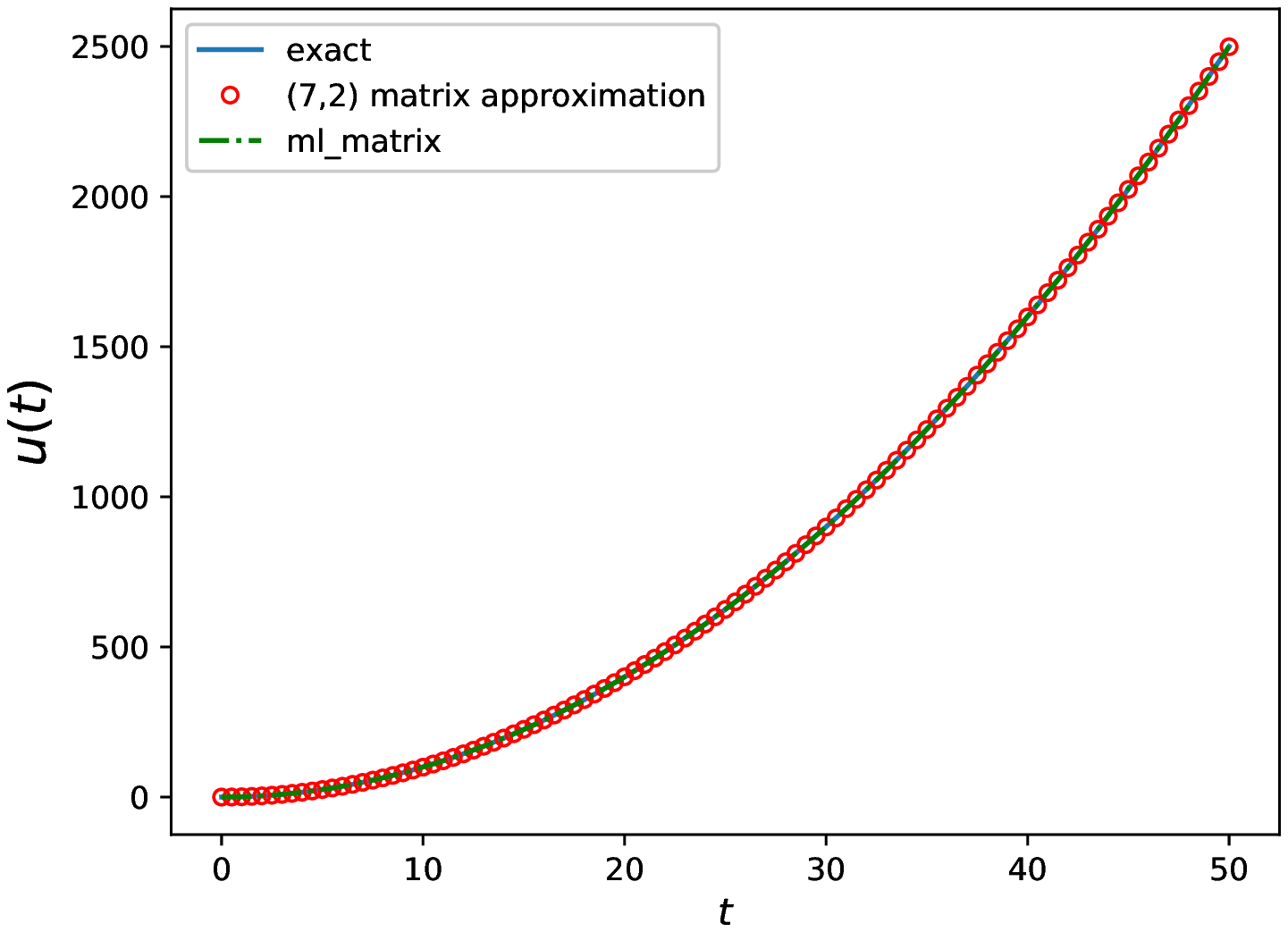}
\caption{
Plots of the exact solution of Bagley-Torvik problem \eqref{eq:Bagley-Torvik} vs the  approximations via $R_{\alpha,\beta}^{7,2}$ and ml\_matrix Matlab function. 
Short-time and long-time profiles are shown in the left and right figures, respectively.
} 
\label{fig:Bagley-Torvik}
\end{figure}

\begin{figure}
\centering
\includegraphics[width=0.49\textwidth]{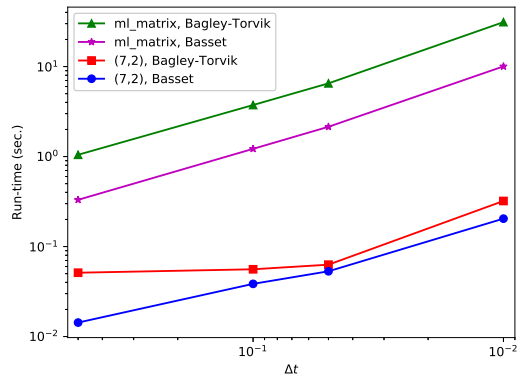}
\caption{
The runtime in seconds for computing the solution of Bagley-Torvik and Basset problems over the interval $[0,50]$ with mesh grid size $\Delta t = 0.01, 0.05, 0.1, 0.5$.
}
\label{fig:cpu-comparison2}
\end{figure}

Figure \ref{fig:Bagley-Torvik} contains a comparison of the profiles of $u$ and its approximation when $a_1 = 3$ and $a_2 = 1$. 
As can be observed, the partial fraction decomposition of $R_{\alpha, \beta}^{7,2}$ generates an accurate approximation of the exact solution.
Furthermore, it follows from  and 
Table \ref{tab:Bagley-Torvik comparison} and
Figure \ref{fig:cpu-comparison2}
that these approximations are effective in terms of runtime.

\begin{table}
\centering
\begin{tabular}{c| cc|c }
	& AE  & RE &  Runtime
	\\ \hline 
	& & &  \\
	$R_{\alpha,\beta}^{7,2}$ 
	& 1.01 & 3.40e-03 &  0.32       
	\\ & & &\\
	matrix ml 
	& 0.05 & 4.83e-05 & 31.16 
	\\ 
\end{tabular}
\caption{
	The maximum absolute error (AE), maximum relative error (RE), and runtime for computing the solution of Bagley-Torvik problem 
	\eqref{eq:Bagley-Torvik} over the interval $[0,50]$ with mesh grid size of $0.01$.
} 	
\label{tab:Bagley-Torvik comparison}	
\end{table}

\subsubsection{The Basset problem}
\label{sec:Basset Problem}
The Basset problem considered in \cite{Mainardi:1995}
\begin{equation}
\label{eq:Basset}
\left[D + \delta^{1-\alpha}\, \OcD^\alpha + 1\right]u(t) = 1, \qquad u(0) = 0, \quad \delta > 0, \quad \alpha \in (0,1).
\end{equation} 
This problem models the dynamics of a sphere immersed in an incompressible viscous
fluid subjected to gravity under a hydraulic force.
When $\alpha$ is a rational number, $\alpha = p/q$, the exact solution to this problem is 
\begin{equation}\label{eq:Basset solution}
u(t) = 1 - \sum\limits_{k=1}^q c_k E_{\frac{1}{q}}(a_kt^\frac{1}{q}),
\end{equation} 
where
$\{a_k\}$ are the zeros of the polynomial 
$x^q + \delta^{(1-\frac{p}{q})}x^p + 1$, 
$c_k = -A_k/a_k$, and 
$A_k = 1/ \prod\limits_{j=1}^q(a_k - a_j)$, $j \ne k$.

When $\alpha = 1/2$, problem \eqref{eq:Basset} can be converted into the system 
\begin{equation}
\label{eq:Basset system}
\OcD^{\frac{1}{2}} U(t) = A U(t) + e_2, \quad U(0) = U^0,
\end{equation} 
where
$$
U(t) = (u, \OcD^\frac{1}{2} u)^T, \qquad  U^0  =(0, 0)^T, \qquad e_2 = (0, 1)^T,
$$ 
and
$$ 
A = \begin{bmatrix}
0  & 1       \\ 
-1 & -\delta^{1/2}      
\end{bmatrix}.
$$
The exact solution of this system (see \cite{Popolizio:2018}) is
\begin{equation}
\label{eq:Basset system solution}
U(t) = t^\frac{1}{2}E_{\frac{1}{2},\frac{3}{2}}(At^\frac{1}{2}) \, e_2,
\end{equation}
which could be approximated by the partial fraction decomposition of $R^{7,2}_{\alpha, \beta}$.

\begin{figure}
\includegraphics[width=0.49\textwidth]{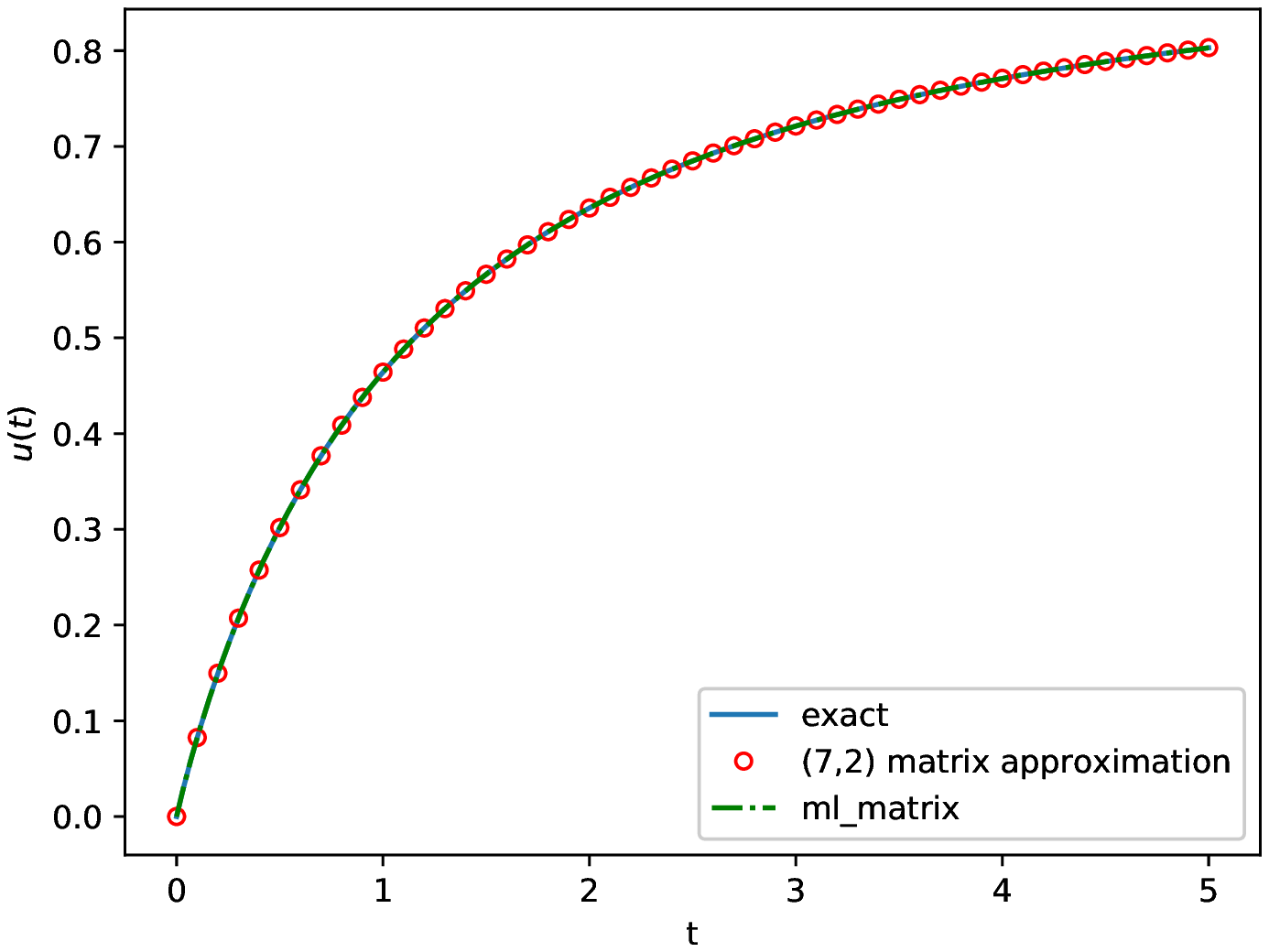}
\includegraphics[width=0.49\textwidth]{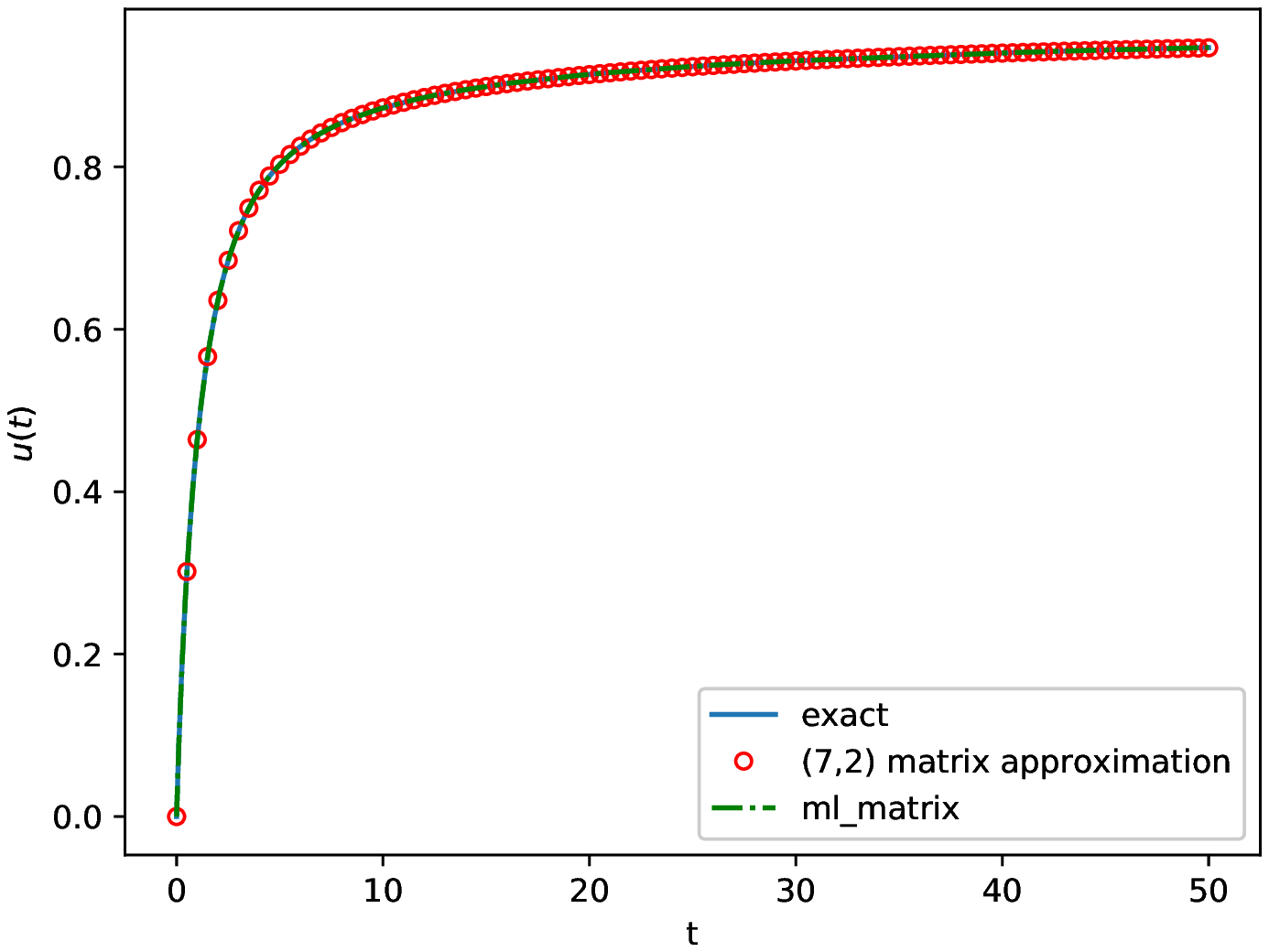}
\caption{
Plots of the exact solution of Basset problem \eqref{eq:Basset solution} vs the  approximations via $R_{\alpha,\beta}^{7,2}$ and ml\_matrix Matlab function. 
Short-time and long-time profiles are shown in the left and right figures, respectively.
}
\label{fig:Basset}
\end{figure}

A comparison between the exact solution \eqref{eq:Basset solution} and the approximation of 
\eqref{eq:Basset system solution} for $\delta = 3/7$ is presented in Figure \ref{fig:Basset}.
With a maximum relative error of 3.17e-04, over the time interval [0,50] for $\Delta t = 0.01$, the runtime of the 
$R^{7,2}_{\alpha,\beta}$ approximation is 0.20 seconds whereas the runtime of the matlab function ml\_matrix is 10.04 seconds. 
A detailed comparison of the runtime is included in Figure~\ref{fig:cpu-comparison2}.

\section{Concluding remarks}

\begin{itemize}

\item
A unified framework for constructing global Pad\'e approximants of the two-parametric MLF, $E_{\alpha, \beta}(-x)$, $x > 0$, 
is presented. In particular, fourth-order global Pad\'e approximants are constructed and tested.

\item
The numerical experiments indicate that these approximants provide efficient and accurate formulas for computing MLFs.
Furhtermore, these approximants perform well when used to approximate the MLF of a matrix.

\item
An algorithm for approximating the inverse function based on the inversion of fourth-order approximants is provided.

\item
The developed rational approximants will play a pivotal role in developing efficient high-order generalized ETD schemes analogue 
to the ETD schemes in \cite{Furati:2018}. This will be explored in a future work.

\end{itemize}

\end{document}